\documentclass[a4paper]{article}%
\usepackage{amsmath}
\usepackage{amsfonts}
\usepackage{amssymb}
\usepackage{graphicx}
\usepackage{enumitem}
\usepackage[hidelinks]{hyperref}
\providecommand{\U}[1]{\protect\rule{.1in}{.1in}}
\newtheorem{theorem}{Theorem}

\newtheorem{definition}{Definition}

\newtheorem{lemma}{Lemma}

\newtheorem{remark}{Remark}

\newenvironment{proof}[1][Proof]{\noindent\textbf{#1.} }{\ \rule{0.5em}{0.5em}}
\newcommand{\norm}[1]{\left\|#1\right\|}
\newcommand{\bra}[1]{\left\langle#1\right\rangle}

\begin{document}
	
	\title{Space-time almost periodicity for the defocusing NLS hierarchy
		via local Toeplitz operators}
	\author{Shuo Zhang$^{1}$}
	\date{}
	\maketitle
	
	\begin{abstract}
		We study the defocusing nonlinear Schr\"odinger hierarchy
		for reflectionless initial data with Dirac spectrum $E$.
		Under a gap moment condition of order $m$ and the assumption that $E$
		is homogeneous, we prove existence and space-time almost periodicity
		of a global classical solution through order $m$ for every such initial
		datum.
		More generally, it suffices that $\mathbb C\setminus E$ be a Dirichlet
		regular Widom domain satisfying the direct Cauchy theorem and that the
		gap maxima of a weighted Green potential be summable, in addition to
		the gap moment condition.
		We obtain the frequencies by summing the upper and lower boundary
		values of the corrected powers $z^r$ on each gap and use local Toeplitz
		operators to reconstruct the potential from the gap phases.
		We also prove uniqueness under bounds on the spatial derivatives.
	\end{abstract}
	
	\section{Introduction}
	\footnotetext[1]{School of Math. Nanjing Univ. Nanjing China, shuozhang@smail.nju.edu.cn}
	We consider the defocusing nonlinear Schr\"odinger hierarchy on the real
	line with reflectionless initial data. In our normalization $t_1=x$ and
	the $t_2$-equation is
	\[
	i\partial_{t_2}q=-\dfrac12\partial_x^2q+|q|^2q.
	\]
	We give conditions on the spectrum $E$ of the associated Dirac operator
	under which the first $m$ positive flows have a common solution
	\[
	q(\mathbf t,x)=F_{q_0}\!\left(x\alpha+\sum_{r=2}^m t_r\beta^{(r)}\right),
	\qquad \mathbf t=(t_2,\ldots,t_m).
	\]
	Here $F_{q_0}$ is continuous on a product of circles, one for each gap of
	$E$, and the frequencies depend only on $E$. In particular, $q$ is a
	uniform limit on $\mathbb R^m$ of trigonometric polynomials in $(x,\mathbf t)$.
	
	For KdV, Deift \cite{Deift08,Deift17} asked whether almost periodic
	initial data give rise to a global solution which is almost periodic in
	time. The periodic theory of McKean and Trubowitz \cite{MT} and
	finite-gap integration \cite{BBEIM} suggest such a conclusion, since the
	evolution is described by linear motion on a torus. With infinitely many
	gaps one must still prove that the potential depends continuously on the
	phases in a topology that gives uniform almost periodicity.
	
	For NLS, Fillman, Li, Luki\'c, and Zhou \cite{FLLZ} proved existence,
	uniqueness, and space-time almost periodicity for uniformly almost periodic
	initial data whose Dirac spectrum is purely absolutely continuous and
	satisfies Craig-type conditions. They linearize the spatial and NLS flows
	by a generalized Abel map. An additional circle in their parametrization
	records rotation of the complex phase. Their direct spectral analysis
	verifies the hypotheses for small analytic quasiperiodic initial data
	with Diophantine frequencies.
	
	In joint work with Kotani, Wang, and Xu \cite{KWXZ}, the author constructed
	global NLS flows for classes containing sufficiently smooth bounded initial
	data, without a decay assumption. The Toeplitz reconstruction in that work
	also gives preservation of spatial Bohr almost periodicity. Here we attach a local Toeplitz operator to each gap of $E$ and use the
	phases of the resulting jump matrices to represent spatial shifts and
	the hierarchy flows as translations on a common torus. This representation
	yields almost periodicity in time for each flow.
	
	On a homogeneous spectrum it suffices to assume a finite moment of the
	gap lengths. We first state the more general criterion in terms of the
	Green function of $\mathbb C\setminus E$.
	
	Fix an index set $\mathcal{J} \subset \mathbb{Z}_{+}$. For $j \in \mathcal{J}$, let
	\[
	I_{j}=[l_{j},r_{j}]=[c_{j}-a_{j},c_{j}+a_{j}], \qquad a_{j}>0, \quad c_{j} \in \mathbb{R},
	\]
	and
	\[
	G = \bigcup_{j \in \mathcal{J}} I_{j}^{\circ}, \qquad E=\mathbb{R} \setminus G, \qquad \Omega=\mathbb{C} \setminus E.
	\]
	Assume that the intervals $I_j$ are pairwise disjoint. The domain assumptions below
	concern the case $G\ne\varnothing$. If $G=\varnothing$, the phase torus is a
	point and the only reflectionless potential is zero \cite[Theorem~6.7]{CG}.
	For $G\ne\varnothing$, let $g_\Omega$ be the positive Green function of $\Omega$, normalized by
	\[
	-\Delta_{z} g_{\Omega}(z,\zeta)=2\pi \delta_{\zeta}, \qquad g_{\Omega}(z,\zeta)=-\log |z-\zeta|+O(1) \quad (z \to \zeta).
	\]
	We call $\Omega$ \emph{Dirichlet regular} if
	\[
	\lim_{z \to e,z \in \Omega} g_{\Omega}(z,\zeta)=0
	\]
	for every $\zeta\in\Omega$ and every boundary point $e$ in the Riemann sphere. We call $\Omega$ a \emph{Widom domain} if
	\[
	\sum_{j} \sup_{s \in I_{j}^{\circ}} g_{\Omega}(s,z)< \infty \quad (z \in \mathbb{C} \setminus \mathbb{R}).
	\]
	A single-valued holomorphic function $F$ on $\Omega$ belongs to the Smirnov class $\mathcal{N}_{+}(\Omega)$ if
	\[
	F \circ \pi = \dfrac{g}{h} , \qquad g,h \in H^{\infty}(\mathbb{D}),
	\]
	where $\pi$ is the universal holomorphic covering map from the unit disk onto $\Omega$ and $h$ is an outer function.
	
	\begin{definition}\label{def:dct}
		Fix $\xi_{*} \in \mathbb{R} \setminus E$. We say that $\Omega$ satisfies the direct Cauchy theorem {\rm(DCT)} if every $F\in\mathcal{N}_{+}(\Omega)$ satisfying
		\[
		\int_{E} \dfrac{|F_{+}(s)|+|F_{-}(s)|}{|s-\xi_{*}|^{2}} ds<\infty
		\]
		obeys
		\[
		F(z)-F(\xi_{*}) = \dfrac{1}{2\pi i}\int_E
		\bigl(F_+(s)-F_-(s)\bigr)
		\left(\frac{1}{s-z}-\frac{1}{s-\xi_*}\right)ds
		\quad(z\in\Omega).
		\]
		Here $F_{\pm}(s)=F(s\pm i0)$ denote the nontangential boundary values. This condition is independent of the choice of the reference point $\xi_{*} \in \mathbb{R} \setminus E$.
	\end{definition}
	
	Define the Green integral by
	\[
	(\mathcal{K}_{\Omega}f)(z):=\dfrac{1}{\pi} \int_{G} g_{\Omega}(z,t)f(t)dt, \qquad z\in \Omega.
	\]
	Set $\omega(s)=(1+s^{2})^{-1/2}$, $v_{\omega}=\mathcal{K}_{\Omega}(\omega|_{G})$ and
	\begin{equation}
		\mathcal{P}_{1}(E) = \sum_{j} \sup_{s \in I_{j}^{\circ}} v_{\omega}(s). \label{1.4}
	\end{equation}
	
	Write
	\begin{equation}
		\lambda_j=1+|c_j|,\quad S_p=\sum_{j\in\mathcal J}a_j\lambda_j^p,
		\quad \Theta=\prod_{j\in\mathcal J}\mathbb T,\quad
		\mathbb T=\mathbb R/(2\pi\mathbb Z). \label{1.1}
	\end{equation}
	Equip $\Theta$ with the translation-invariant metric
	\[
	d_{\Theta}(\theta,\varphi)=\sum_{j \in \mathcal{J}} 2^{-j}\left|e^{i\theta_{j}}-e^{i\varphi_{j}}\right|.
	\]
	It induces the product topology, and $\Theta$ is compact.
	
	For $q\in L^1_{\mathrm{loc}}(\mathbb R;\mathbb C)$, let $D_q$ be the
	self-adjoint Dirac operator in $L^2(\mathbb R;\mathbb C^2)$ associated with
	\[
	D_q=iJ\frac{d}{dx}+
	\begin{pmatrix}0&q\\\overline q&0\end{pmatrix},
	\qquad J=\operatorname{diag}(1,-1),
	\]
	on its maximal domain. In the Schur normalization, the Weyl functions are given by
	$m_-^q(z)=f_{-,2}(0,z)/f_{-,1}(0,z)$ for $z\in\mathbb C_+$ and
	$m_+^q(z)=f_{+,2}(0,z)/f_{+,1}(0,z)$ for $z\in\mathbb C_-$, where
	$f_\pm(\cdot,z)$ are nonzero solutions of $D_qf=zf$ square integrable at
	$\pm\infty$, respectively. At a general base point $x$, we write
	$m_\pm^q(x,z)$ for the corresponding ratios at $x$. In $\mathbb C_+$
	we also use the right Weyl coordinate
	$m_+^q(x,z)=\overline{m_+^q(x,\bar z)}=f_{+,1}(x,z)/f_{+,2}(x,z)$.
	Define
	\[
	\begin{split}
		\mathcal R(E)=\bigl\{q\in L^1_{\mathrm{loc}}(\mathbb R;\mathbb C):\ &
		\sigma(D_q)=E,\\
		&m_-^q(\lambda+i0)=m_+^q(\lambda-i0)
		\text{ for a.e. }\lambda\in E\bigr\}.
	\end{split}
	\]
	
	We fix the normalization of the defocusing AKNS hierarchy
	(cf.~\cite[Appendix~B]{FLLZ}). For smooth $q$, put
	\[
	U=-izJ+P(q),\qquad
	P(q)=\begin{pmatrix}0&-i\overline q\\ iq&0\end{pmatrix}.
	\]
	The coefficients of the formal series $R=J+\sum_{r\ge1}R_rz^{-r}$ are
	determined successively by $R_x=[R,U]$ and $R^2=I_2$, where
	$[A,B]=AB-BA$. The $r$th flow is $\partial_{t_r}q=\mathcal X_r(q)$, with
	\[
	\mathcal X_r(q)=\frac1i[J,R_{r+1}]_{21},\qquad r\ge0 .
	\]
	The degree-zero flow is $q\mapsto e^{-2it_0}q$. Each $\mathcal X_r$ is
	a differential polynomial of order $r$. In particular,
	\[
	\mathcal X_1(q)=q_x,\quad
	\mathcal X_2(q)=\dfrac i2 q_{xx}-i|q|^2q,\quad
	\mathcal X_3(q)=-\dfrac14 q_{xxx}+\dfrac32|q|^2q_x.
	\]
	
	\begin{theorem}\label{t1}
		Let $m\geq 2$ be an integer. Suppose that $\Omega$ is a Dirichlet regular Widom domain with DCT and that 
		\[
		S_{m}<\infty, \qquad \mathcal{P}_{1}(E)<\infty.
		\]
		There exist real sequences $\beta^{(r)}\in\mathbb R^{\mathcal J}$,
		$1\le r\le m$, depending only on $E$, such that
		\begin{equation}
			|\beta_j^{(r)}|\le C_r\lambda_j^r,\qquad
			\alpha=\beta^{(1)},\quad \beta=\beta^{(2)}. \label{1.2}
		\end{equation}
		For each $q_0\in\mathcal R(E)$ there is a continuous function
		$F_{q_0}:\Theta\to\mathbb C$ such that $x\mapsto F_{q_0}(x\alpha)$ is a
		$C^m$ representative of $q_0$. With this representative,
		\[
		q(\mathbf t,x)=F_{q_0}\!\left(x\alpha+\sum_{r=2}^m t_r\beta^{(r)}\right)
		\]
		is a global classical solution of
		\begin{equation}
			\partial_{t_r}q=\mathcal X_r(q)\quad (2\le r\le m),
			\qquad q(\mathbf0,x)=q_0(x). \label{1.3}
		\end{equation}
		The phases are taken modulo $2\pi$. All derivatives
		\[
		\partial_x^k\partial_{t_2}^{\ell_2}\cdots\partial_{t_m}^{\ell_m}q,
		\qquad k+\sum_{r=2}^m r\ell_r\le m,
		\]
		are jointly Bohr almost periodic in $(x,\mathbf t)$.
		The solution of \eqref{1.3} is unique among classical
		solutions $v$ with
		\[
		\sup_{\substack{|\mathbf t|\leq T\\x\in\mathbb R}}
		\sum_{k=0}^r|\partial_x^k v(\mathbf t,x)|<\infty\qquad(T>0), \quad (k=0,1 \text{ for the cubic flow}).
		\]
	\end{theorem}
	
	\begin{theorem}\label{t2}
		Suppose that $G \neq \varnothing$, $S_{0}<\infty$ and $E$ is homogeneous, namely, there exists a constant $\eta>0$ such that
		\[
		|E \cap (x-r,x+r)| \geq \eta r \quad (\text{for all } x \in E,\ r>0).
		\]
		Then $\Omega$ is a Dirichlet regular Widom domain with DCT and
		\[
		\mathcal{P}_{1}(E) \leq C_{\eta}S_{0}.
		\]
	\end{theorem}
	
	For homogeneous $E$, the two theorems leave just one summability
	condition, $S_m<\infty$. Every $q_0\in\mathcal R(E)$ then has a $C^m$
	almost periodic representative, and the same gap moment gives joint
	almost periodicity of all mixed derivatives of weighted order
	$k+\sum_{r=2}^m r\ell_r\le m$.
	
	The frequencies can be read off before any initial potential is chosen.
	For $1\le r\le m$, the sum of the two boundary values of $\Phi_{r,N}$
	on the $j$th retained gap $I_{j}^{\circ}$ is a constant $\beta^{(r)}_{j,N}$. It is
	obtained by correcting $z^r$ with the Cauchy integral of the Green
	potential of $rs^{r-1}$ on the retained gaps (Lemma~\ref{l3}). As
	$N\to\infty$, $\beta^{(r)}_{j,N}\to\beta^{(r)}_j$. The Weyl functions
	of $q_0$ determine the remaining part of the jump. The $r$th flow thus
	adds $t_r\beta_j^{(r)}$ to the phase on the $j$th gap.
	
	For a finite set of gaps the local Toeplitz equations give a
	reconstruction function $F_N$ on their product of circles. The passage
	to infinitely many gaps rests on uniqueness of the normalized matrix
	with the prescribed jumps. Indeed, if $W_1,W_2$ are two admissible
	matrices, $W_2^{-1}W_1$ extends through the open gaps. Its entries
	belong to the Smirnov class and have equal boundary values on the
	spectrum. DCT and the normalization at infinity give $W_1=W_2$
	(Lemma~\ref{l8}). The bound on $\mathcal P_1(E)$ gives compactness
	of the Green potentials in the finite-gap decompositions. Uniqueness
	then identifies the limit even when the phases vary with $N$, and
	\[
	\lim_{N\to\infty}\sup_{\theta\in\Theta}
	|F_N(\theta)-F_{q_0}(\theta)|=0.
	\]
	The finite-gap trace estimates give the spatial derivatives, and the
	hierarchy equations give the time derivatives. For each flow, uniqueness on finite
	time strips follows by comparing the time transfer matrices on the two half-line
	Weyl directions (Lemma~\ref{l12}).
	
	The hypotheses of Fillman, Li, Luki\'c, and Zhou imply that
	$q_0\in\mathcal R(E)$, that $E$ is homogeneous, and that $S_2<\infty$
	\cite[(1.5)--(1.6), Theorem~4.8, and Lemma~7.12]{FLLZ}.
	Thus Theorems~\ref{t1} and~\ref{t2} extend the existence, uniqueness,
	and space-time almost periodicity result of \cite[Theorem~1.2]{FLLZ}
	for the cubic NLS to a larger class of spectra. With
	$S_m<\infty$, the same conclusions hold for each flow of order $r\le m$.
	
	For KdV and Toda the analogous results concern different spectral
	problems. Eichinger, VandenBoom, and Yuditskii \cite{EVY} obtained global
	classical solutions of the KdV hierarchy, uniformly almost periodic in
	space and time. Their construction on an Abelian covering uses DCT, a
	summability condition on the critical values of the Martin function,
	and a gap moment depending on the flow. Vinnikov and Yuditskii \cite{VY}
	linearized the Toda hierarchy on the character group associated with a
	homogeneous Jacobi spectrum. Theorems~\ref{t1} and~\ref{t2} give a corresponding result for the
	defocusing NLS hierarchy on a Dirac spectrum unbounded in both directions.
	Our proof proceeds by a different method. We use local Toeplitz operators
	on the spectral gaps and identify each flow with a translation of the
	phases in the jump matrices.
	
	We recall briefly some earlier results. Sodin and Yuditskii
	\cite{SY95,SY97} constructed isospectral tori for reflectionless
	Sturm--Liouville and Jacobi operators with homogeneous spectra and
	established spatial almost periodicity, assuming finite total gap
	length in the continuum setting. For KdV, Egorova \cite{Egorova}
	proved time almost periodicity for a class of limit-periodic initial
	data, and Craig \cite{Craig} developed trace formulas and equations for
	Dirichlet data. Damanik and Goldstein \cite{DG16} proved global
	existence and uniqueness for small analytic quasiperiodic KdV data
	with Diophantine frequencies. Binder, Damanik, Goldstein, and Luki\'c
	\cite{BDGL} proved time almost periodicity under spectral hypotheses.
	The estimates in \cite{DG14,DGL16,DGL17} give their application to
	small analytic quasiperiodic data. For Toda, Binder, Damanik,
	Luki\'c, and VandenBoom \cite{BDLV} obtained almost periodicity in
	space and time under Craig-type conditions. The gap estimates of
	Leguil, You, Zhao, and Zhou \cite{LYZZ} give applications to subcritical
	analytic one-frequency Schr\"odinger data.
	
	Earlier inverse-spectral treatments of NLS with nondecaying data
	include \cite{BE,BM}. Bessonov, Luki\'c, and Yuditskii \cite{BLY}
	studied spatial almost periodicity for reflectionless canonical
	systems and Dirac operators on Widom domains with DCT. For the
	defocusing mKdV flow see Li and Luki\'c \cite{LL}. Oh \cite{OhAP,OhLP}
	studied local well-posedness in spaces of almost periodic functions and
	proved global existence for a regular class of limit-periodic defocusing
	NLS data. Local quasiperiodic solutions were studied in
	\cite{DLX,Papenburg}.
	
	Without further assumptions, spatial Bohr almost periodicity can be lost
	under the KdV flow, as the bounded solution of Chapouto, Killip, and
	Vi\c{s}an \cite{CKV} shows. Damanik, Li, and Strungaru \cite{DLS}
	obtained Weyl almost periodicity under the Sodin--Yuditskii hypotheses
	of absolute continuity, finite total gap length, and the Parreau--Widom
	condition.
	
	Section~2 constructs the frequencies from Green integrals. Section~3
	solves the finite systems of Toeplitz equations and proves the trace
	estimates used in the limit. In Section~4 we prove uniqueness of the
	admissible matrix, pass to infinitely many gaps, and prove
	Theorems~\ref{t1} and~\ref{t2}.
	
	\section{Green integrals and frequencies}
	
	\subsection{Green integrals}
	Throughout this section, $\Omega$ is Dirichlet regular and $S_0<\infty$.
	Recall the Green integral
	\[
	(\mathcal{K}_{\Omega}f)(z)=\frac{1}{\pi}\int_G g_{\Omega}(z,t)f(t)\,dt,\qquad z\in\Omega .
	\]
	For a measure $\mu$ on $G$ we use the same notation, with $f(t)\,dt$ replaced by $d\mu(t)$, whenever the integral exists.
	We shall use the symmetry and domain monotonicity of the Green function,
	\begin{equation}
		g_\Omega(z,w)=g_\Omega(w,z)\ge0,\qquad
		\Omega_1\subset\Omega_2\ \Longrightarrow\ g_{\Omega_1}\le g_{\Omega_2}.
		\label{2.15}
	\end{equation}
	Throught the paper, the kernel of the smaller domain is extended by zero if either argument lies outside that domain.
	If $\mu$ is a nonnegative Radon measure on $G$ and $P=\mathcal K_\Omega\mu\not\equiv+\infty$, the Riesz decomposition gives
	\begin{equation}
		h\text{ harmonic in }\Omega,\quad 0\le h\le P
		\quad\Longrightarrow\quad h\equiv0 .
		\label{2.1}
	\end{equation}
	
	For a real function $f$ on $G$, write
	\begin{equation}
		v_f(s)=
		\begin{cases}
			(\mathcal K_\Omega f)(s),&s\in G,\\
			0,&s\in E
		\end{cases}
		\label{2.9}
	\end{equation}
	whenever the integral exists. Let $\mathcal J_N=\{j\in\mathcal J:j\le N\}$,
	$G_N=\bigcup_{j\in\mathcal J_N}I_j^o$, $E_N=\mathbb R\setminus G_N$ and
	$\Omega_N=\mathbb C\setminus E_N$. Write $K_N=\overline{G_N}$ and
	$\partial K_N=\{l_j,r_j:j\in\mathcal J_N\}$ for its endpoints.
	Define $v_{f,N}$ by \eqref{2.9} with $\Omega_N,G_N,E_N$ in place of $\Omega,G,E$.
	
	\begin{lemma}\label{l1}
		For every real $f\in L^\infty(G)$,
		\[
		\sup_{z\in\Omega}|(\mathcal K_\Omega f)(z)|\le M\|f\|_{L^\infty(G)},
		\qquad M=\frac{\pi S_0}{2}.
		\]
		Moreover, $v_f\in C(I_j)$ and $v_f(l_j)=v_f(r_j)=0$ for every $j$. In particular,
		\[
		0\le v_\omega\le v_1\le M,\qquad
		v_{1,N}\uparrow v_1,\qquad v_{\omega,N}\uparrow v_\omega .
		\]
	\end{lemma}
	
	\begin{proof}
		Choose $I=(c-a,c+a)$ containing $\overline{G_N}$, and put
		\[
		u_I(z)=\Im\bigl(\sqrt{(z-c)^2-a^2}-(z-c)\bigr),\qquad z\in\mathbb C_+,
		\]
		where $\sqrt{(z-c)^2-a^2}\sim z-c$ at infinity. Extend $u_I$ by reflection across the real axis. Its boundary values on $I$ are $\sqrt{a^2-(s-c)^2}$, and it vanishes on $\mathbb R\setminus I$ and at infinity. In particular $0\le u_I\le a$. Since $\partial_yu_I(s+i0)=-1$ for $s\in I$,
		$-\Delta u_I=2\boldsymbol1_I\,ds$ in $\Omega_I=\mathbb C\setminus(\mathbb R\setminus I)$, where $ds$ denotes length measure on the real axis.
		In the Riesz decomposition of $u_I$, the nonnegative harmonic part is bounded above by $u_I$ and hence vanishes by the maximum principle. Thus
		\[
		0\le\mathcal K_{\Omega_N}1\le\mathcal K_{\Omega_I}1=u_I\le a .
		\]
		
		Let $0\le f\le A$ and $V_N=\mathcal K_{\Omega_N}f$, so that $0\le V_N\le B_N:=aA$.
		For $e\in E_N$ and $\delta>0$, put $B=B(e,\delta)$, $D_N=B\cap\Omega_N$ and $\Gamma_N=\partial B\cap\Omega_N$.
		The local Green representation, with $\omega_{D_N}^z$ denoting harmonic measure, gives
		\[
		g_{\Omega_N}(z,t)=g_{D_N}(z,t)
		+\int_{\Gamma_N}g_{\Omega_N}(\zeta,t)\,d\omega_{D_N}^z(\zeta),
		\qquad z\in D_N .
		\]
		Integrating in $t$ and using $g_{D_N}\le g_B$ yields
		\begin{equation}
			V_N(z)\le\frac A\pi\int_{e-\delta}^{e+\delta}g_B(z,t)\,dt
			+B_N\omega_{D_N}^z(\Gamma_N)
			\le CA\delta+B_N\omega_{D_N}^z(\Gamma_N).
			\label{2.2}
		\end{equation}
		Here $C=2(1+\log2)/\pi$, since $g_B(z,t)\le\log(2\delta/|z-t|)$. Each finite boundary point of $\Omega_N$ is regular. Thus $\omega_{D_N}^z(\Gamma_N)\to0$ as $z\to e$. Letting $\delta\downarrow0$ in \eqref{2.2} gives $V_N(z)\to0$ at $e$. At an interior point $s_0\in G_N$, choose $\chi\in C_c^\infty(G_N)$ equal to one near $s_0$. In a disk about $s_0$,
		\begin{equation}
			V_N(z)=-\frac1\pi\int_{G_N}\chi(t) f(t)\log|z-t|\,dt+h(z)
			\label{2.3}
		\end{equation}
		with $h$ harmonic. The logarithmic integral is continuous for bounded $f$. This proves the asserted continuity for $v_{f,N}$, also for signed $f$ by taking positive and negative parts.
		
		Assume $G_N\ne\varnothing$ and write $v_N=v_{1,N}$. Its maximum $m_N>0$ is attained at a point $s_N\in G_N$.
		The Poisson formula for bounded harmonic functions \cite[Chapter~7]{ABR} gives
		\begin{equation}
			(\mathcal K_{\Omega_N}f)(x+iy)=(P_y*v_{f,N})(x),
			\qquad P_y(x)=\frac{y}{\pi(x^2+y^2)},\quad y>0 .
			\label{2.5}
		\end{equation}
		Near $s_N$, $\mathcal K_{\Omega_N}1+|y|$ is harmonic and even in $y$. Hence
		$-\partial_y\mathcal K_{\Omega_N}1(s_N+i0)=1$, and
		\begin{align*}
			1&=\lim_{y\downarrow0}\frac{m_N-(\mathcal K_{\Omega_N}1)(s_N+iy)}{y}\\
			&=\lim_{y\downarrow0}\frac1\pi\int_{\mathbb R}
			\frac{m_N-v_N(t)}{(s_N-t)^2+y^2}\,dt
			\ge\frac{m_N}{\pi}\int_{E_N}\frac{dt}{(s_N-t)^2}
			\ge\frac{4m_N}{\pi|G_N|}.
		\end{align*}
		The last integral is minimized, among complements of sets of length $|G_N|$, when the deleted set is the interval of that length centered at $s_N$.
		By \eqref{2.5} and symmetry,
		\[
		0\le\mathcal K_{\Omega_N}1(z)\le\frac{\pi|G_N|}{4}\le M.
		\]
		As $\Omega_N\uparrow\Omega$, the Green kernels increase to $g_\Omega$ \cite[Theorem~5.2]{Pap}. Monotone convergence gives the two asserted limits and
		$|\mathcal K_\Omega f|\le\|f\|_\infty\mathcal K_\Omega1\le M\|f\|_\infty$.
		
		Finally, for $e\in E$, put $D=B\cap\Omega$. We may take $B_N=AM$ in \eqref{2.2}. Domain monotonicity of the exit probability gives
		$\omega_{D_N}^z(\Gamma_N)\le\omega_D^z(\partial B\cap\Omega)$. Passing to the limit in \eqref{2.2} for $0\le f\le A$ gives
		\[
		0\le(\mathcal K_\Omega f)(z)\le CA\delta+AM\omega_D^z(\partial B\cap\Omega).
		\]
		Since $e$ is regular, the last term tends to zero as $z\to e$. Letting $\delta\downarrow0$ gives the endpoint values. The local decomposition \eqref{2.3} applies equally in $\Omega$ and proves continuity in each open gap.
	\end{proof}
	
	\subsection{Frequencies}
	We use the Hilbert and Cauchy transforms $\mathcal{M}$
	\[
	Hv(s)=\frac1\pi\operatorname{pv}\int_{\mathbb R}\frac{v(t)}{s-t}\,dt,
	\qquad Cv(z)=\frac1\pi\int_{\mathbb R}\frac{v(t)}{t-z}\,dt,\quad z\in\mathbb C\setminus\mathbb R .
	\]
	For $v\in L^2(\mathbb R)$ their boundary values satisfy
	\begin{equation}
		(Cv)(s\pm i0)=-Hv(s)\pm i v(s)\qquad\text{for a.e. }s\in\mathbb R .
		\label{2.10}
	\end{equation}
	Write $\mathcal D'(G)$ for the space of distributions on $G$.
	
	\begin{lemma}\label{l2}
		For every real $f\in L^1(G)\cap L^2(G)$,
		\begin{equation}
			\|v_f\|_{L^p(\mathbb R)}\le M\|f\|_{L^p(G)},\qquad
			v_{f,N}\longrightarrow v_f\quad\text{in }L^p(\mathbb R),\quad p=1,2,
			\label{2.4}
		\end{equation}
		and
		\begin{equation}
			\partial_s Hv_f=f \quad\text{in }\mathcal D'(G).
			\label{2.8}
		\end{equation}
		Thus, if $p$ is locally absolutely continuous on each $I_j^o$ and $p'=f$ there, then
		$p-Hv_f=b_j$ a.e. on $I_j^o$ for a real constant $b_j$.
	\end{lemma}
	
	\begin{proof}
		By symmetry and Lemma~\ref{l1}, the integral of $g_\Omega(s,t)/\pi$ over either real variable, which is $v_{1}$, is bounded by  $M$. Tonelli's theorem and the Cauchy--Schwarz inequality give
		\[
		\|v_f\|_1\le M\|f\|_1,\qquad
		\|v_f\|_2^2\le M\int_G |f(t)|^2v_1(t)\,dt\le M^2\|f\|_2^2.
		\]
		By monotone convergence of the extended kernels, for a.e. $s\in G$,
		\[
		|v_f(s)-v_{f,N}(s)|\le\frac1\pi\int_G(g_\Omega-g_{\Omega_N})(s,t)|f(t)|\,dt
		\longrightarrow0.
		\]
		The right-hand side is bounded by $v_{|f|}\in L^1\cap L^2$, so dominated convergence proves \eqref{2.4}.
		
		For bounded $f$, Lemma~\ref{l1} and the Poisson representation give
		\begin{equation}
			(\mathcal K_\Omega f)(x+iy)=(P_y*v_f)(x),\qquad y>0 .
			\label{2.6}
		\end{equation}
		For general $f$, truncate its positive and negative parts and apply monotone convergence. Since $\int_{\mathbb R}(P_y*v_{|f|})(x)\,dx=\|v_{|f|}\|_1$, the potential is locally area integrable. Fubini's theorem then gives $-\Delta\mathcal K_\Omega f=2f\,ds$ in $\Omega$.
		
		Fix $\varphi\in C_c^\infty(G)$, and let $\chi\in C_c^\infty(G)$ equal one near $\operatorname{supp}\varphi$.
		The local decomposition \eqref{2.3} holds near $\operatorname{supp}\varphi$, with $h$ harmonic and even in $y$. There
		\[
		-\partial_y(\mathcal K_\Omega f)(x+iy)
		=(P_y*(\chi f))(x)-\partial_y h(x+iy).
		\]
		Since $P_y*(\chi f)\to\chi f$ in $L^1$ and $\partial_yh(x+iy)\to0$ uniformly on $\operatorname{supp}\varphi$,
		\begin{equation}
			\lim_{y\downarrow0}\int_{\mathbb R}
			(-\partial_y\mathcal K_\Omega f)(x+iy)\varphi(x)\,dx
			=\int_G f(s)\varphi(s)\,ds .
			\label{2.7}
		\end{equation}
		The Fourier multipliers of $H$ and $P_y*$ are $-i\operatorname{sgn}\xi$ and $e^{-y|\xi|}$, respectively. Hence
		\[
		-\partial_y(P_y*v_f)=\partial_x(P_y*Hv_f),
		\qquad P_y*Hv_f\to Hv_f\quad\text{in }L^2 .
		\]
		Combining \eqref{2.6} and \eqref{2.7}, we obtain
		$\int_G f\varphi=-\int_{\mathbb R}(Hv_f)\varphi'$, which is \eqref{2.8}. \bigskip
	\end{proof}
	
	\begin{lemma}\label{l3}
		Suppose $S_m<\infty$, $m\ge2$. For $1\le r\le m$, let $p_r(s)=s^r$, $f_r=p_r'$, $v_r=v_{f_r}$ and $v_{r,N}=v_{f_r,N}$. Then $v_{r,N}\to v_r$ in $L^1(\mathbb R)$. Put
		\[
		\Phi_r(z)=p_r(z)+\frac1\pi\int_G\frac{v_r(s)}{s-z}\,ds,\qquad
		\Phi_{r,N}(z)=p_r(z)+\frac1\pi\int_{G_N}\frac{v_{r,N}(s)}{s-z}\,ds.
		\]
		There are real numbers $\beta_j^{(r)}$ such that
		\begin{equation}
			\Phi_r(s+i0)+\Phi_r(s-i0)=\beta_j^{(r)},\qquad s\in I_j^o .
			\label{2.13}
		\end{equation}
		The same relation defines $\beta_{j,N}^{(r)}$ for $j\in\mathcal J_N$, with $\beta_{j,N}^{(r)}=0$ for $j\notin\mathcal J_N$. We write $\beta^{(1)}=\alpha$, $\beta^{(2)}=\beta$ and likewise for the truncations.
		For each fixed $j$, $\beta_{j,N}^{(r)}\to\beta_j^{(r)}$.
		Moreover, $\Phi_{r,N}\to\Phi_r$ locally uniformly on $\mathbb C\setminus\mathbb R$, and the frequencies satisfy \eqref{1.2}.
	\end{lemma}
	
	\begin{proof}
		Since $|s|\le(1+S_0)\lambda_j$ on $I_j$,
		\[
		\|f_r\|_{L^1(G)}\le 2r(1+S_0)^{r-1}S_{r-1}.
		\]
		The $L^1$ estimate and dominated convergence in the proof of Lemma~\ref{l2} use only $f_r\in L^1(G)$, so $\|v_{r,N}-v_r\|_1\to 0$. On each $G_N$, $f_r$ is bounded, and \eqref{2.8} gives $\partial_s(p_r-Hv_{r,N})=0$ in $\mathcal D'(G_N)$. For $\varphi\in C_c^\infty(G)$,
		\[
		\int_{\mathbb R} (Hv_{r,N})\varphi
		=-\int_{\mathbb R} v_{r,N}H\varphi
		\longrightarrow-\int_{\mathbb R}v_rH\varphi
		=\langle Hv_r,\varphi\rangle ,
		\]
		since $H\varphi$ is bounded. Thus $\partial_s(p_r-Hv_r)=0$ in $\mathcal D'(G)$.
		The Poisson representation \eqref{2.6} and $-\Delta\mathcal K_\Omega f_r=2f_r\,ds$ follow by truncation in $L^1$ as in Lemma~\ref{l2}. Near each interior point of a gap,
		\[
		(\mathcal K_\Omega f_r)(x+iy)+\operatorname{sgn}(y)\Im p_r(x+iy)
		\]
		is harmonic, since the second term has distributional Laplacian $2p'_r(x)\,ds$. Thus $v_r$ is smooth there. Cutting off the Cauchy integral near the point gives $\Phi_r(s\pm i0)=p_r(s)-Hv_r(s)\pm iv_r(s)$ pointwise. The same argument applies to $\Phi_{r,N}$, proving \eqref{2.13} and its finite analogue.
		
		Fix $j$ and choose $\varphi_j\in C_c^\infty(I_j^o)$ with $\int\varphi_j=1$. For $N\ge j$,
		\[
		|\beta_{j,N}^{(r)}-\beta_j^{(r)}|
		=2\left|\int_{\mathbb R}(v_{r,N}-v_r)H\varphi_j\right|
		\le2\|v_{r,N}-v_r\|_1\|H\varphi_j\|_\infty\longrightarrow0 .
		\]
		For compact $K\subset\mathbb C\setminus\mathbb R$,
		\[
		\sup_{z\in K}|\Phi_{r,N}(z)-\Phi_r(z)|
		\le\frac{\|v_{r,N}-v_r\|_1}{\pi\operatorname{dist}(K,\mathbb R)}
		\longrightarrow0 .
		\]
		
		Fix $N$ and let $D_R=\Omega_N\cap B(0,R)$, where $G_N\Subset(-R,R)$.
		For $t\in G_N$ the Green representation gives
		\[
		g_{D_R}(s,t)=\int_{\partial D_R}\log|s-\zeta|\,d\omega_{D_R}^t(\zeta)-\log|s-t|.
		\]
		This identity extends to the zero extension of $g_{D_R}(\cdot,t)$ on the real axis, by harmonicity of $\log|s-\cdot|$ outside its pole and continuity at the regular boundary.
		Taking the Hilbert transform in $s$ gives, almost everywhere,
		\[
		H\!\left(\frac{g_{D_R}(\cdot,t)}{\pi}\right)(s)=
		\frac12\operatorname{sgn}(s-t)
		-\frac1\pi\int_{\partial D_R}
		\arctan\frac{s-\Re\zeta}{|\Im\zeta|}\,d\omega_{D_R}^t(\zeta).
		\]
		Here the odd conjugate of $\log|s-\zeta|$ about $s=\Re\zeta$ is $-\arctan((s-\Re\zeta)/|\Im\zeta|)$. At $\Im\zeta=0$ the arctangent is interpreted as $(\pi/2)\operatorname{sgn}(s-\zeta)$. The formula holds first in distributions, since the two logarithmic terms have equal mass. In particular,
		\begin{equation}
			\left|H\!\left(\frac{g_{D_R}(\cdot,t)}{\pi}\right)(s)\right|\le1 .
			\label{2.11}
		\end{equation}
		For $f=f_r$, put $v_{f,R}(s)=\pi^{-1}\int_{G_N}g_{D_R}(s,t)f(t)\,dt$.
		Integrating \eqref{2.11} against $|f(t)|\,dt$ gives $\|Hv_{f,R}\|_\infty\le\|f\|_{L^1(G_N)}$. As $R\to\infty$, domain monotonicity and the bound by $v_{|f|,N}$ give $v_{f,R}\to v_{f,N}$ in $L^2$. Then one can choose an almost everywhere convergent subsequence, passing to the limit gives
		\begin{equation}
			\|Hv_{f,N}\|_\infty\le\int_{G_N}|f(t)|\,dt .
			\label{2.12}
		\end{equation}
		Evaluating the constant gap average at $s=c_j$ gives
		\begin{equation}
			|\beta_{j,N}^{(r)}-2c_j^r|\le 2\int_{G_N}|f_r(t)|\,dt,\qquad j\in\mathcal J_N .
			\label{2.16}
		\end{equation}
		Letting $N\to\infty$, we obtain
		\begin{equation}
			|\beta_j^{(r)}|\le 2|c_j|^r+4r(1+S_0)^{r-1}S_{r-1}\le C_r\lambda_j^r .
			\label{2.18}
		\end{equation}
		\bigskip
	\end{proof}
	
	Along the imaginary axis,
	\[
	|\Phi_r(iy)-p_r(iy)|\le\frac{\|v_r\|_1}{\pi|y|},\qquad y\ne0 .
	\]
	For fixed $N$, $v_{r,N}$ has compact support, so
	\[
	\Phi_{r,N}(z)=p_r(z)-\frac1{\pi z}\int_{\mathbb R}v_{r,N}(s)\,ds+O(z^{-2}),
	\qquad |z|\to\infty .
	\]
	Moreover, Lemma~\ref{l1} and \eqref{2.12} give
	\[
	\|v_{r,N}\|_\infty\le M\|f_r\|_{L^\infty(G_N)},\qquad
	\|Hv_{r,N}\|_\infty\le\|f_r\|_{L^1(G_N)} .
	\]
	Since $\Phi_{r,N}(x+iy)-p_r(x+iy)=-(P_y*Hv_{r,N})(x)+i(P_y*v_{r,N})(x)$ for $y>0$, the difference is bounded on $\mathbb C_+$, and also on $\mathbb C_-$ by reflection. In particular, for each $e\in\partial K_N$,
	\begin{equation}
		\Phi_{r,N}(z)=O(1)\qquad(z\to e).
	\end{equation}
	
	\section{Finite truncation and associated results}
	Throughout the paper, for a complex function $f$ and $s \in \mathbb{R}$, we denote
	\[
	f_{+}(s)=f(s+i0)=\lim_{y \downarrow 0}f(s+iy), \qquad f_{-}(s)=f(s-i0)=\lim_{y \downarrow 0}f(s-iy).
	\]
	
	\subsection{The initial Weyl function}
	Fix $q_0\in\mathcal R(E)$ and put
	\begin{equation}
		m(z)=\begin{cases}
			m_-^{q_0}(z),&z\in\mathbb C_+,\\
			m_+^{q_0}(z),&z\in\mathbb C_-,
		\end{cases}
		\qquad m^\#(z)=\overline{m(\overline z)}. \label{3.24}
	\end{equation}
	The Weyl functions in this normalization satisfy $|m(z)|<1$ and
	$m(iy)\to0$ as $y\to\pm\infty$ by \cite[Sections~2--3]{CG}.
	The reflectionless identity glues these bounded holomorphic functions
	across every open interval in $E^\circ$. Indeed, their boundary values
	agree almost everywhere, so Morera's theorem, applied across the real
	axis, gives a holomorphic extension. The maximum principle then yields
	$|m|<1$ on $\mathbb C\setminus\overline G$.
	
	Set
	\[
	p=mm^\#,\qquad \Delta=1-p,\qquad
	M=\begin{pmatrix}1&m\\m^\#&1\end{pmatrix},\qquad
	\mathcal{M}_{q_{0}}=\Delta^{-1/2}M.
	\]
	Since $|p|<1$, we have $\operatorname{Re}\Delta>0$ and use the square
	root with positive real part. Thus
	\[
	\det \mathcal{M}_{q_{0}}=1,\qquad \mathcal{M}_{q_{0}}(iy)\to I_2\quad(y\to\pm\infty),\qquad
	\frac{(\mathcal{M}_{q_{0}})_{12}}{(\mathcal{M}_{q_{0}})_{11}}=m.
	\]
	
	In a full-line gap $I_j^\circ$, the half-line Weyl functions are
	meromorphic. Hence, for either
	half-line realization, each $\lambda\in I_j^\circ$ is a resolvent point
	or an isolated eigenvalue. In both cases $D_q f=\lambda f$ has a nonzero
	solution square integrable at that end. Since $f^*Jf=0$, we may normalize
	the two Weyl solutions as
	\[
	f_L(0,\lambda)=(1,e^{iu_j(\lambda)})^T,\qquad
	f_R(0,\lambda)=(1,e^{iv_j(\lambda)})^T,
	\]
	where $f_L\in L^2(\mathbb R_-;\mathbb C^2)$ and
	$f_R\in L^2(\mathbb R_+;\mathbb C^2)$.
	We therefore have
	\begin{equation}
		m(\lambda+i0)=e^{iu_j(\lambda)},\qquad
		m(\lambda-i0)=e^{iv_j(\lambda)}, \quad \text{on } I_j^\circ. \label{3.25}
	\end{equation}
	The Lagrange identity
	\[
	\partial_x(f^*J\partial_\lambda f)=-i\norm f^2
	\]
	gives
	\[
	u_j'(\lambda)=\int_{-\infty}^0\norm{f_L(x,\lambda)}^2\,dx>0,
	\qquad
	v_j'(\lambda)=-\int_0^\infty\norm{f_R(x,\lambda)}^2\,dx<0.
	\]
	The two Weyl directions cannot coincide in a spectral gap.
	We may consequently choose their continuous arguments so that
	\[
	0<u_j-v_j<2\pi\quad\text{on }I_j^\circ.
	\]
	The monotone functions $u_j$ and $v_j$ have finite one-sided limits at both endpoints, since the total increase of $u_j-v_j$ is at most $2\pi$. Extend them to $I_j$ by these limits. For $\phi_j=(u_j+v_j)/2$, we have
	\begin{equation}
		\operatorname{Var}_{I_{j}} \phi_{j} \leq \dfrac{1}{2} \int_{I_j^\circ}(u_{j}'-v_{j}')ds \leq \pi.
		\label{3.1}
	\end{equation}
	Writing $d=u_j-v_j\in(0,2\pi)$, we obtain
	\[
	\Delta_+=1-e^{id},\qquad
	\Delta_-=\overline{\Delta_+},\qquad
	\arg\Delta_+=\frac d2-\frac\pi2,
	\]
	and
	\[
	M_+=\begin{pmatrix}0&e^{iu_j}\\e^{-iv_j}&0\end{pmatrix}M_-.
	\]
	Hence
	\begin{equation}
		(\mathcal{M}_{q_{0}})_+=iQ_{\phi_j}(\mathcal{M}_{q_{0}})_-
		\quad\text{on }I_j^\circ,\qquad
		Q_\phi=\begin{pmatrix}0&e^{i\phi}\\e^{-i\phi}&0\end{pmatrix}.
		\label{3.2}
	\end{equation}
	
	For $h\in C([0,1])$, define the Bernstein polynomial
	\[
	B_nh(t)=\sum_{k=0}^n h(k/n)\binom nk t^k(1-t)^{n-k}.
	\]
	It satisfies $B_nh(0)=h(0)$ and $B_nh(1)=h(1)$, and
	\[
	(B_nh)'(t)=n\sum_{k=0}^{n-1}\left(h\left(\frac{k+1}n\right)-h\left(\frac kn\right)\right)\binom{n-1}{k} t^k(1-t)^{n-1-k}.
	\]
	Thus $B_n$ preserves monotonicity and the total variation of monotone functions.
	
	Let
	\[
	U_{j}(t) = u_{j}(l_{j}+2a_{j}t), \qquad V_{j}(t) = v_{j}(l_{j}+2a_{j}t), \qquad t_{j}(z) = \dfrac{z-l_{j}}{2a_{j}},
	\]
	and
	\[
	\begin{array}{c}
		u_{j,N}(z) = B_{N}U_{j}(t_{j}(z)), \qquad v_{j,N}(z) = B_{N}V_{j}(t_{j}(z)), \\
		\phi_{j,N}(z) = \dfrac{u_{j,N}(z)+v_{j,N}(z)}{2}.
	\end{array}
	\]
	By monotonicity and preservation of endpoint values,
	\[
	\operatorname{Var}_{I_j}\phi_{j,N}
	\leq\dfrac12\bigl(\operatorname{Var}_{I_j} u_j+\operatorname{Var}_{I_j} v_j\bigr)\leq\pi.
	\]
	By the classical Bernstein approximation theorem \cite[Theorem~1]{VV2016}, we have
	\[
	\phi_{j,N}\longrightarrow\phi_j\quad\text{in }C(I_j).
	\]
	On a compact subinterval of $I_j^\circ$, choose a smooth cutoff equal to one on a neighborhood of that subinterval. The derivative approximation theorem \cite[Theorem~2]{VV2016} applies to the cutoff functions, and the discarded Bernstein weights and their derivatives of any fixed order are exponentially small on the subinterval. Thus
	\begin{equation}
		\phi_{j,N}\longrightarrow\phi_j \quad\text{in }C^\infty_{\mathrm{loc}}(I_j^\circ). \label{3.28}
	\end{equation}
	
	\subsection{Local Toeplitz operators}
	For $j\in\mathcal J_N$, let $s_j(z)=\sqrt{(z-l_j)(z-r_j)}$, with $s_j(z)\sim z-c_j$ at infinity, and write
	\[
	\bra{f}_{j}:=\dfrac{1}{\pi} \int_{I_{j}} \dfrac{f(s)}{\rho_{j}(s)}ds , \qquad \rho_{j}(s) = \sqrt{a_{j}^{2}-(s-c_{j})^{2}}.
	\]
	Set
	\[
	\vartheta_{j,N}=\bra{\phi_{j,N}}_j,\qquad \psi_{j,N}=\phi_{j,N}-\vartheta_{j,N},\qquad \bra{\psi_{j,N}}_j=0,
	\]
	and
	\[
	\xi_{j,N}(z)=\frac{s_j(z)}{2}\bra{\frac{\psi_{j,N}(s)}{z-s}}_j.
	\]
	Then
	\[
	\xi_{j,N,+}+\xi_{j,N,-} = \psi_{j,N}, \qquad \xi_{j,N}(z) = O_{j}(z^{-1}) \quad (z \to \infty).
	\]
	Since $\bra{1}_{j}=1$ and $\bra{(z-s)^{-1}}_{j}=s_{j}(z)^{-1}$, in a neighborhood of $I_j$,
	\[
	\xi_{j,N}(z)=\frac{\psi_{j,N}(z)}2-
	\frac{s_j(z)}{2\pi}\int_{I_j}
	\frac{\psi_{j,N}(z)-\psi_{j,N}(\lambda)}{z-\lambda}
	\frac{d\lambda}{\rho_j(\lambda)}.
	\]
	Thus, for $e\in\partial I_j$,
	\[
	\xi_{j,N}(e)=\frac{\psi_{j,N}(e)}2,\qquad
	\xi_{j,N}(z)-\xi_{j,N}(e)=O(|z-e|^{1/2}).
	\]
	
	Set
	\[
	\chi_j(z)=\left(\frac{z-c_j-a_j}{z-c_j+a_j}\right)^{1/4},
	\qquad \chi_j(z)\to1\quad(z\to\infty),
	\]
	and
	\[
	\widehat H_j(\phi,z)=\frac{\chi_j+\chi_j^{-1}}2I_2
	+\frac{\chi_j-\chi_j^{-1}}2Q_\phi.
	\]
	Then $\det\widehat H_j=1$ and
	$\widehat H_{j,+}=iQ_\phi\widehat H_{j,-}$.
	Define
	\[
	E_{j,N}(z)=e^{i\xi_{j,N}(z)J},\qquad
	H_{j,N}(\theta_j,z)=E_{j,N}(z)\widehat H_j(\vartheta_{j,N}+\theta_j,z).
	\]
	Using $e^{iwJ}Q_\phi=Q_\phi e^{-iwJ}$, we obtain
	\[
	\det H_{j,N}=1,\qquad
	(H_{j,N})_+=iQ_{\phi_{j,N}+\theta_j}(H_{j,N})_-,\qquad
	H_{j,N}(z)=I_2+O_j(z^{-1}).
	\]
	Moreover, $H_{j,N},H_{j,N}^{-1} = O(|z-e|^{-1/4})$.
	
	For $j\in\mathcal J_N$, choose $d_j>a_j$ so that the disks $D_j=\{z:|z-c_j|<d_j\}$ have pairwise disjoint closures. Orient $\Gamma_j=\partial D_j$ counterclockwise. Let $\mathcal A(\overline D_j)$ be the Banach space of functions continuous on $\overline D_j$ and holomorphic in $D_j$, with the supremum norm. Let $P_j^+$ and $P_j^-$ be the interior and exterior Cauchy projections
	on $\Gamma_j$, with the convention
	\[
	\begin{aligned}
		P_j^+v(z)&=\frac1{2\pi i}\int_{\Gamma_j}
		\frac{v(\zeta)}{\zeta-z}\,d\zeta,
		&&z\in D_j,\\
		P_j^-v(z)&=\frac1{2\pi i}\int_{\Gamma_j}
		\frac{v(\zeta)}{z-\zeta}\,d\zeta,
		&&z\notin\overline D_j.
	\end{aligned}
	\]
	Their boundary values satisfy $P_j^++P_j^-=I$ whenever these are defined. Define
	\[
	T_j(H_{j,N})f=P_j^+(H_{j,N}f), \qquad \omega_{j,N}=I_{2}-H_{j,N}.
	\]
	To see that this is bounded on $\mathcal A(\overline D_j)$, put $\Gamma_j'=\{|z-c_j|=(a_{j}+d_{j})/2\}$.
	Contour deformation gives
	\[
	T_j(H_{j,N})f(z)=f(z)+\frac1{2\pi i}\int_{\Gamma_j'}
	(H_{j,N}(\zeta)-I_2)\frac{f(z)-f(\zeta)}{z-\zeta}\,d\zeta.
	\]
	The quotient has a removable singularity at $z=\zeta$.
	For $z\in\Gamma_j$, 
	\[
	\norm{T_j(H_{j,N})f(z)}\leq C \norm f_{\infty,D_j}.
	\]
	The maximum principle gives the same bound throughout $D_j$. The argument also applies to $H_{j,N}^{-1}$.
	
	For each $\theta\in\Theta$, we seek a solution of
	\begin{equation}
		T_i(H_{i,N})f_{i,N}=I_2-\sum_{j \in \mathcal{J}_{N},j\ne i}P_j^-(\omega_{j,N}f_{j,N}), \quad i \in \mathcal{J}_{N},
		\label{3.3}
	\end{equation}
	in $X_{N}:=\prod_{j \in \mathcal{J}_{N}} \mathcal A(\overline D_j;\mathbb{C}^{2 \times 2})$ with norm
	\[
	\norm{f}_{X_{N}}=\max_{j \in \mathcal{J}_{N}} \sup_{z \in \overline{D}_{j}} \norm{f_{j}(z)}.
	\]
	Since $H_{i,N},H_{i,N}^{-1}$ are holomorphic off $I_i$ and tend to $I_2$
	at infinity, exterior factors give the Toeplitz identities
	\[
	T_i(H_{i,N})^{-1}=T_i(H_{i,N}^{-1}),\qquad
	T_i(H_{i,N}^{-1})I_2=I_2.
	\]
	For $j\ne i$, the Cauchy formula also gives
	\[
	T_i(H_{i,N}^{-1})P_j^-(\omega_{j,N}f_{j,N})
	=P_j^-(H_{i,N}^{-1}\omega_{j,N}f_{j,N})\quad\text{on }D_i.
	\]
	Thus \eqref{3.3} is equivalent to $f=\boldsymbol I_N+\mathcal K_{N,\theta} f$, where
	\begin{equation}
		(\mathcal K_{N,\theta} f)_i(z)=\frac1{2\pi i}\sum_{j \in \mathcal{J}_{N},j\ne i}
		\int_{\Gamma_j}
		\frac{H_{i,N}(\theta_i,\zeta)^{-1}\omega_{j,N}(\theta_j,\zeta)f_j(\zeta)}
		{\zeta-z}\,d\zeta, \label{3.4}
	\end{equation}
	for $z\in\mathbb C\setminus\overline{\bigcup_{j\ne i}D_j}$ and $i \in \mathcal{J}_{N}$.
	
	\begin{lemma}
		\label{l4}For every $N$ and $\theta\in\Theta$, $I-\mathcal K_{N,\theta}$ is boundedly invertible on $X_N$. In particular,
		\[
		f_{N,\theta} = (I-\mathcal{K}_{N,\theta})^{-1}\boldsymbol{I}_{N}, \qquad \boldsymbol{I}_{N} = (I_{2})_{j \in \mathcal{J}_{N}}
		\]
		is the unique solution to \eqref{3.3} in $X_{N}$.
	\end{lemma}
	
	\begin{proof}
		Fix $N,\theta$. All indices in the proof belong to $\mathcal J_N$, and we suppress the subscript $N$ on the local factors and phases.
		
		Since the $D_j$ are separated from one another, the right-hand side of \eqref{3.4} is holomorphic in a neighborhood $U_{i}$ of $\overline{D}_{i}$ for every $i$, and
		\[
		\sup_{z\in U_{i}} \norm{(\mathcal{K}_{\theta}f)_{i}(z)} \leq C_{\theta} \norm{f}_{X_{N}}.
		\]
		Cauchy estimates and the Arzel\`a--Ascoli theorem yield the compactness of $\mathcal{K}_{\theta}$. Thus by the Fredholm alternative, it suffices to show $\ker (I-\mathcal{K}_{\theta}) = \{0\}$.
		
		Let $f\in\ker(I-\mathcal K_\theta)$. Then
		\begin{equation}
			T_{i}(H_{i})f_{i} = - \sum_{j \neq i} P_{j}^{-} (\omega_{j}f_{j}), \label{3.5}
		\end{equation}
		and we set
		\[
		Z(z) = \left\{
		\begin{array}{ll}
			-\sum_{j}P_{j}^{-}(\omega_{j}f_{j})(z), & z \notin \bigcup_{j} \overline{D}_{j}, \\
			H_{i}(z)f_{i}(z), & z \in D_{i} \setminus I_{i}.
		\end{array}
		\right.
		\]
		By \eqref{3.5} and $P_{j}^{+}+P_{j}^{-}=I$, $Z$ extends holomorphically across all $\Gamma_{j}$, is holomorphic on $\mathbb{C} \setminus \bigcup_{j}I_{j}$ and satisfies
		\[
		Z_{+} = VZ_{-} \quad a.e. \text{ on } \mathbb{R}, \qquad V(s) = \left\{
		\begin{array}{ll}
			iQ_{\phi_{j}(s)+\theta_{j}}, & s \in I_{j}^\circ, \\
			I_{2}, & s \in E_{N}.
		\end{array}
		\right.
		\]
		Furthermore,
		\[
		Z(z) = O(z^{-1}) \quad (z \to \infty), \qquad Z(z)=O(|z-e|^{-1/4}) \quad (z \to e \in \bigcup_{j}\partial I_{j}).
		\]
		
		The function $B(z)=Z(\bar z)^*Z(z)$ is holomorphic in $\mathbb C_+$. For large $R$ and small $\tau>0$, apply Cauchy's theorem on
		\[
		U_{\tau,R} = \{ z:\Im z >0,|z|<R,|z-e|>\tau \text{ for every } e \in \partial K_{N} \}.
		\]
		Since $B(z)=O(|z-e|^{-1/2})$ at each endpoint and $B(z)=O(z^{-2})$ at infinity, the integrals over the small semicircles and the large semicircle are $O(\sqrt\tau)$ and $O(R^{-1})$, respectively. Letting $\tau\downarrow0$ and $R\to\infty$ gives
		\[
		\int_{\mathbb{R}} Z_{-}^{*}Z_{+} ds=0.
		\]
		Taking Hermitian parts and using $V+V^*=2\boldsymbol{1}_{E_N}I_2$, we obtain
		\[
		\int_{E_{N}}Z_{-}^{*}Z_{-}ds = \int_{\mathbb{R}} Z_{-}^{*} \dfrac{V+V^{*}}{2} Z_{-}ds =0.
		\]
		Thus $Z$ vanishes almost everywhere on $E_N$. Since it is holomorphic across $E_N^\circ$, the identity theorem gives $Z\equiv0$ on $\mathbb C\setminus K_N$. The invertibility of each $H_i$ gives $f=0$. \bigskip
	\end{proof}
	
	\subsection{Reconstruction and the NLS hierarchy}
	Denote the unique solution of \eqref{3.3} by $f_{N,\theta} = (f_{j,N,\theta})_{j \in \mathcal{J}_{N}}$. Define
	\[
	W_{N}(\theta,z)=\begin{cases}
		H_{i,N}(\theta_i,z)f_{i,N,\theta}(z),&z\in D_i\setminus I_i,\\
		I_2-\displaystyle\sum_jP_j^-(\omega_{j,N}f_{j,N,\theta})(z),
		&z\notin\displaystyle\bigcup_j\overline D_j.
	\end{cases}
	\]
	Then $W_{N}$ is holomorphic on $\mathbb C\setminus \bigcup_{j \in \mathcal{J}_{N}}I_{j}$ and
	\begin{equation}
		W_{N,+}(\theta,\lambda)=iQ_{\phi_{j,N}(\lambda)+\theta_j}
		W_{N,-}(\theta,\lambda)\quad(\lambda\in I_j^\circ).
		\label{3.6}
	\end{equation}
	Since $\det(iQ_\phi)=1$, $\det W_N$ has no jump on $I_{j}^{\circ}$, is $O(|z-e|^{-1/2})$ at each endpoint and tends to $1$ at infinity. Its endpoint singularities are removable, so Liouville's theorem and the adjugate formula give
	\begin{equation}
		\begin{array}{c}
			\det W_{N} \equiv 1, \qquad W_{N}(\theta,z)=I_{2}+\dfrac{A_{N}(\theta)}{z}+O(z^{-2}), \\
			W_{N},W_{N}^{-1}=O(|z-e|^{-1/4}) \quad (e \in \partial K_{N}).
		\end{array} \label{3.7}
	\end{equation}
	Expanding the exterior Cauchy integral gives $A_N$. We set
	\begin{equation}
		A_{N}(\theta) = -\dfrac{1}{2\pi i} \sum_{j} \int_{\Gamma_{j}} \omega_{j,N}(\zeta)f_{j,N,\theta}(\zeta) d\zeta, \qquad F_{N}(\theta)=2(A_{N}(\theta))_{21}. \label{3.27}
	\end{equation}
	Among matrices holomorphic on $\mathbb{C} \setminus \overline{G_{N}}$
	that extend continuously to each open gap from each of the two
	half-planes, \eqref{3.6} and \eqref{3.7} determine $W_N$ uniquely. Indeed, for any other such solution $\widetilde W_N$, the quotient $W_N^{-1}\widetilde W_N$ has no jump on $I_{j}^{\circ}$ and has growth $O(|z-e|^{-1/2})$ at the endpoints. It is entire and tends to $I_2$ at infinity. Applying uniqueness to $Q_0\overline{W_N(\theta,\bar z)}Q_0$ gives
	\begin{equation}
		W_{N}(\theta,z) = Q_{0}\overline{W_{N}(\theta,\bar{z})}Q_{0}, \qquad (A_{N})_{12}=\overline{(A_{N})_{21}}, \quad (A_{N})_{22} = \overline{(A_{N})_{11}}. \label{3.8}
	\end{equation}
	The map $\theta\mapsto\mathcal K_{N,\theta}$ is continuous in operator norm. By continuity of inversion and compactness of the finite torus,
	\begin{equation}
		\sup_{\theta \in \mathbb{T}^{\mathcal{J}_{N}}} \norm{(I-\mathcal{K}_{N,\theta})^{-1}}_{X_{N}\to X_{N}}<\infty. \label{3.9}
	\end{equation}
	\begin{remark}
		\label{r1} The factors in \eqref{3.4} depend smoothly on the retained phases. The identity
		\[
		\partial_{\theta} (I-\mathcal{K})^{-1}=(I-\mathcal{K})^{-1}(\partial_{\theta} \mathcal{K}) (I-\mathcal{K})^{-1}
		\]
		shows that $f_{N,\theta}$, $A_N(\theta)$ and $F_N(\theta)$ are smooth functions of $(\theta_j)_{j\in\mathcal J_N}$ and $2\pi$-periodic in each coordinate.
	\end{remark}
	
	\begin{lemma}
		\label{l5}For $\mathbf t=(t_2,\ldots,t_m)$, let
		\[
		h_N=x\alpha_N+\sum_{r=2}^m t_r\beta_N^{(r)},\qquad
		q_{N,\theta}(\mathbf t,x)=F_N(\theta+h_N).
		\]
		Then $\partial_{t_r}q_{N,\theta}=\mathcal X_r(q_{N,\theta})$ for $2\le r\le m$. In particular,
		\[
		i\partial_{t_2}q=-\dfrac{1}{2}\partial_{x}^{2}q+|q|^{2}q.
		\]
	\end{lemma}
	
	\begin{proof}
		Write $W=W_N(\theta+h_N,z)$ and set
		\[
		\Phi=x\Phi_{1,N}(z)+\sum_{r=2}^m t_r\Phi_{r,N}(z),\qquad
		\Psi=e^{-i\Phi J}W .
		\]
		Suppress $N$ on the frequencies and phases below. On $I_j^\circ$, \eqref{2.13} and \eqref{3.6} give $\Phi_++\Phi_-=h_j$ and
		\begin{equation}
			\Psi_+ =e^{-i\Phi_+J}e^{ih_jJ}iQ_{\phi_j+\theta_{j}}W_- =e^{i\Phi_-J}iQ_{\phi_j+\theta_{j}}W_-
			=iQ_{\phi_j+\theta_{j}}\Psi_-. \label{3.10}
		\end{equation}
		Thus the jump of $\Psi$ is independent of $x,\mathbf t$. Put $R=W^{-1}JW$ and define
		\begin{equation}
			\begin{aligned}
				U&=\Psi^{-1}\partial_x\Psi
				=W^{-1}\partial_xW-i\Phi_{1}R,\\
				V_r&=\Psi^{-1}\partial_{t_r}\Psi
				=W^{-1}\partial_{t_r}W-i\Phi_{r}R,\qquad 2\le r\le m.
			\end{aligned} \label{3.11}
		\end{equation}
		By \eqref{3.10}, $U$ and $V_r$ extend across the open gaps $I_{j}^{\circ}$. In $D_j$, the factorization $W=H_jf_j$ and smooth dependence on the phases give $W_x,W_{t_r}=O(|z-e|^{-1/4})$, just as for $W,W^{-1}$. Since each $\Phi_r$ is bounded at the endpoints, $U,V_r=O(|z-e|^{-1/2})$ there. Thus $U$ and $V_r$ are entire.
		
		Write $W=I_2+A/z+B/z^2+O(z^{-3})$ and $R=\sum_{n\ge0}R_n z^{-n}$, $R_0=J$. Since $\Phi_r=z^r+O(z^{-1})$, comparison of the polynomial parts in \eqref{3.11} gives
		\begin{equation}
			\begin{aligned}
				U&=-izJ+P,\qquad P=-i[J,A]=\begin{pmatrix}0&-i\overline q\\iq&0\end{pmatrix},\\
				V_r&=-i(z^rR)_{\ge0}=-i\sum_{n=0}^r z^{r-n}R_n .
			\end{aligned} \label{3.12}
		\end{equation}
		From $R_x=[R,U]$ and $R^2=I_2$ we have, for $n\ge0$,
		\[
		i[J,R_{n+1}]=(R_n)_x-[R_n,P],\qquad
		JR_{n+1}+R_{n+1}J=-\sum_{a=1}^n R_aR_{n+1-a}.
		\]
		The first identity determines the off-diagonal part of $R_{n+1}$, the second its diagonal part. Thus these are the normalized hierarchy coefficients and $R_{r+1}$ is a differential polynomial in $q,\bar q$ of order at most $r$. The constant term of $\partial_{t_r}U-\partial_xV_r+[V_r,U]=0$ gives
		\[
		P_{t_r}=-i\bigl((R_r)_x-[R_r,P]\bigr)=[J,R_{r+1}],
		\qquad \partial_{t_r}q=\mathcal X_r(q).
		\]
		For $r=2$, \eqref{3.12} gives $V_2=zU-iR_2$. Using $R_1=[J,A]$, $R_2=[J,B]-A[J,A]$ and the coefficient of $z^{-1}$ in $W_x=WU+i\Phi_1JW$, we obtain
		\[
		-iR_2=-i[J,B]-AP=-A_x+i\kappa_1J,\qquad
		\Phi_1=z+\kappa_1/z+O(z^{-2}).
		\]
		The diagonal entries of $AP$ are $i|q|^2/2$ and $-i|q|^2/2$. The off-diagonal entries of $-iR_2$ are determined by $A_{21}=q/2$ and $A_{12}=\overline q/2$. Thus
		\[
		V_2=zU+\begin{pmatrix}
			-i|q|^2/2&-\partial_x\overline q/2\\
			-\partial_xq/2&i|q|^2/2
		\end{pmatrix}.
		\]
		The $(2,1)$ entry of the zero-curvature equation is
		\[
		iq_{t_2}-\left(izq_x-\tfrac12q_{xx}\right)
		+\left(izq_x-|q|^2q\right)=0,
		\]
		which is the stated NLS equation. \bigskip
	\end{proof}
	
	Fix $N,\theta$ and write $W_{N,\theta}(x,z)=W_N(\theta+x\alpha_N,z)$ and $q_{N,\theta}(x)=F_N(\theta+x\alpha_N)$. At $\mathbf t=0$,
	\[
	\partial_{x} \Psi_{N,\theta}(x,z) = \Psi_{N,\theta}(x,z) U_{N,\theta}(x,z), \qquad U_{N,\theta}^{T} = \begin{pmatrix}
		-iz & iq \\
		-i\overline{q} & iz
	\end{pmatrix}.
	\]
	Thus, if
	\[
	Y_{N,\theta}(x,z) = (\boldsymbol{\psi}_{L,N,\theta}(x,z),\boldsymbol{\psi}_{R,N,\theta}(x,z)) = \Psi_{N,\theta}(x,z)^{T},
	\]
	then its columns solve $D_{q_{N,\theta}}\boldsymbol\psi=z\boldsymbol\psi$ and are given by
	\[
	\begin{array}{c}
		\boldsymbol{\psi}_{L,N,\theta}(x,z) = e^{-ix\Phi_{1,N}(z)}\begin{pmatrix}
			W_{N,\theta,11}(x,z) \\
			W_{N,\theta,12}(x,z)
		\end{pmatrix}, \\
		\boldsymbol{\psi}_{R,N,\theta}(x,z) = e^{ix\Phi_{1,N}(z)}\begin{pmatrix}
			W_{N,\theta,21}(x,z) \\
			W_{N,\theta,22}(x,z)
		\end{pmatrix}.
	\end{array}
	\]
	For each $z\in\mathbb C_+$, compactness of the finite torus gives
	\[
	\sup_{x \in \mathbb{R}} \left( \norm{W_{N,\theta}(x,z)}+\norm{W_{N,\theta}(x,z)^{-1}} \right)<\infty.
	\]
	Thus,
	\[
	\norm{\boldsymbol{\psi}_{L,N,\theta}(x,z)} \leq Ce^{\gamma x} \quad (x \leq 0), \qquad \norm{\boldsymbol{\psi}_{R,N,\theta}(x,z)} \leq Ce^{-\gamma x} \quad (x \geq 0),
	\]
	where $\gamma=\Im\Phi_{1,N}(z)=\Im z+(\mathcal K_{\Omega_N}1)(z)>0$. Therefore,
	\[
	\boldsymbol{\psi}_{L,N,\theta} \in L^{2}(\mathbb{R}_{-}), \qquad \boldsymbol{\psi}_{R,N,\theta} \in L^{2}(\mathbb{R}_{+}).
	\]
	Hence
	\[
	m_{-}^{q_{N,\theta}}(z) = \dfrac{W_{N,\theta,12}(0,z)}{W_{N,\theta,11}(0,z)}, \qquad m_{+}^{q_{N,\theta}}(z) = \dfrac{W_{N,\theta,21}(0,z)}{W_{N,\theta,22}(0,z)}.
	\]
	In the second formula we use the right Weyl coordinate in $\mathbb C_+$, the conjugate reflection of its lower half-plane coordinate. By \eqref{3.8},
	\[
	m_{-}^{q_{N,\theta}}(\lambda+i0)=\overline{m_{+}^{q_{N,\theta}}(\lambda+i0)}, \quad a.e. \ \lambda \in E_{N}.
	\]
	Write $q=q_{N,\theta}(\mathbf 0,\cdot)$ and
	$W(z)=W_{N,\theta}(0,z)$. Let $g_q$ denote half the trace of the
	full-line Weyl matrix. Its representing measure has support
	$\sigma(D_q)$; see \cite[(3.19)]{FLLZ}. In $\mathbb C_+$, substituting $m_-^q=W_{12}/W_{11}$ and $m_+^q=W_{21}/W_{22}$
	and using $\det W=1$, we obtain
	\[
	g_q=\frac{i}{2}\frac{1+m_-^q m_+^q}{1-m_-^q m_+^q}
	=\frac{i}{2}S,\qquad
	S=W_{11}W_{22}+W_{12}W_{21}.
	\]
	On $E_N^\circ$, \eqref{3.8} gives
	$|W_{11}|^2-|W_{12}|^2=1$, so $|m_-^q|<1$. By the reflectionless
	identity,
	\[
	\Im g_q(s+i0)
	=\frac{1+|m_-^q(s+i0)|^2}{2(1-|m_-^q(s+i0)|^2)}>0
	\qquad (s\in E_N^\circ).
	\]
	Thus $E_N=\overline{E_N^\circ}\subset\sigma(D_q)$.
	
	On each $I_{j}^{\circ}$, the jump relation \eqref{3.6} gives $S_+=-S_-$.
	Since $g_q=-iS/2$ in $\mathbb C_-$ by \eqref{3.8} and Schwarz
	symmetry, the two boundary values of $g_q$ agree there. Hence $g_q$
	extends holomorphically across $G_N$ with real values on $G_N$.
	Its representing measure therefore vanishes on $G_N$, which gives
	$G_N\subset\rho(D_q)$. We conclude that $\sigma(D_q)=E_N$.
	Replacing $\theta$ by $\theta+\sum_{r=2}^m t_r\beta_N^{(r)}$ yields
	\begin{equation}
		q_{N,\theta}(\mathbf t,\cdot)\in\mathcal R(E_N)
		\qquad (\mathbf t\in\mathbb R^{m-1}).
		\label{3.22}
	\end{equation}
	
	\subsection{Canonical decomposition and uniform estimates}
	Let $d\mu(s)=\rho(s)ds$ with $\rho \in L^{p}(G_{N},\mathbb{R})$, $1<p<\infty$. If $U \in C(\mathbb{C};\mathbb{R})$ satisfies
	\[
	-\Delta U=2\mu \quad \text{in } \mathcal{D}'(\Omega_{N}), \qquad U|_{E_{N}}=0, \qquad \lim_{z \to \infty} U(z)=0,
	\]
	then
	\begin{equation}
		U(z) = (\mathcal{K}_{\Omega_{N}}\mu)(z) = \dfrac{1}{\pi} \int_{G_{N}} g_{\Omega_{N}}(z,t)d\mu(t), \quad z\in\Omega_{N}. \label{3.17}
	\end{equation}
	Indeed, put $P=\mathcal K_{\Omega_N}\mu$ and $h=U-P$.
	Comparison with the Green function of a single larger gap $I \supset G_{N}$ gives
	$P(z)\to0$ at infinity and, by H\"older's inequality near the logarithmic
	singularity in \eqref{2.2}, at each point of $E_N$. Thus
	\[
	-\Delta P=2\mu,\qquad \Delta h=0,\qquad h|_{E_N}=0,\qquad h(\infty)=0.
	\]
	The maximum principle gives $h=0$, hence \eqref{3.17}.
	
	Write $q_{N,\theta}=q_{N,\theta}(\mathbf 0,\cdot)$ and $W_{N,\theta}(z)=W_N(\theta,z)$. By \eqref{3.22}, $q_{N,\theta}\in\mathcal R(E_N)$. For $z\in\mathbb C_+$, let
	\[
	\Delta_{N,\theta}(z)=1-m_{-}^{q_{N,\theta}}(z)m_{+}^{q_{N,\theta}}(z), \qquad \mathcal{M}_{q_{N,\theta}} = \Delta_{N,\theta}^{-1/2}\begin{pmatrix}
		1 & m_{-}^{q_{N,\theta}} \\
		m_{+}^{q_{N,\theta}} & 1
	\end{pmatrix}.
	\]
	Use the positive-real-part square root and extend $\mathcal{M}_{q_{N,\theta}}$ to $\mathbb C_-$ by $\mathcal{M}_{q_{N,\theta}}(z)=Q_0\overline{\mathcal{M}_{q_{N,\theta}}(\bar z)}Q_0$. Let $\phi_{j,N,\theta}$ be the continuous phase given by \eqref{3.2} for this potential, so that
	\begin{equation}
		\mathcal{M}_{q_{N,\theta},+} = iQ_{\phi_{j,N,\theta}}\mathcal{M}_{q_{N,\theta},-} \quad \text{on } I_j^\circ. \label{3.14}
	\end{equation}
	As in \eqref{3.1}, $\operatorname{Var}_{I_j}\phi_{j,N,\theta}\leq\pi$.
	
	\begin{lemma}
		\label{l6} There is a real-symmetric function $\eta_{N,\theta}$, holomorphic on $\mathbb C\setminus \overline{G_{N}}$ and normalized by $\eta_{N,\theta}(\infty)=0$, such that
		\[
		W_{N,\theta} = e^{i\eta_{N,\theta}J}\mathcal{M}_{q_{N,\theta}}.
		\]
		Put $u_{N,\theta}=\Im(\eta_{N,\theta})_+$ on $G_N$ and $u_{N,\theta}=0$ on $E_N$. There is a real signed measure $\mu_{N,\theta}$ on $G_N$ such that
		\begin{equation}
			\mu_{N,\theta}|_{I_j^\circ} = -\dfrac{1}{2} d(\phi_{j,N}-\phi_{j,N,\theta}), \qquad |\mu_{N,\theta}|(I_{j}) \leq \pi, \label{3.15}
		\end{equation}
		and
		\[
		u_{N,\theta}=(\mathcal{K}_{\Omega_{N}}\mu_{N,\theta})|_{\mathbb{R}}, \qquad \eta_{N,\theta}(z) = \dfrac{1}{\pi} \int_{G_{N}} \dfrac{u_{N,\theta}(s)}{s-z}ds.
		\]
	\end{lemma}
	
	\begin{proof}
		In this proof, write $W=W_{N,\theta}$, $q=q_{N,\theta}$, $\mathcal{M}=\mathcal{M}_q$ and $\Delta=\Delta_{N,\theta}$. Since $\det W=1$,
		\[
		W=\operatorname{diag}(d,d^{-1})\mathcal{M}, \qquad d=W_{11}\sqrt{\Delta}, \quad d^{-1}=W_{22}\sqrt{\Delta},
		\]
		so $d$ and $d^{-1}$ are holomorphic on $\mathbb C_+$ and tend to $1$ at infinity. Choose $\eta=-i\log d$ there with $\eta(\infty)=0$. On $E_N^\circ$, the reflection symmetry of $W,\mathcal{M}$ gives $\overline d=d^{-1}$. Thus $\eta$ is real there, and Schwarz reflection gives a single-valued real-symmetric continuation to $\mathbb C\setminus K_N$.
		
		The Herglotz function
		\[
		g_{q}(z) = \dfrac{i}{2} \dfrac{1+m_{-}^{q}m_{+}^{q}}{1-m_{-}^{q}m_{+}^{q}}
		\]
		is real analytic on $G_N$ and purely imaginary on $E_N^\circ$. The exponential representation of $g_q$ gives
		\[
		g_{q}(z) = \dfrac{i}{2} \prod_{j \in \mathcal{J}_{N}} \dfrac{z-\nu_{j}}{\sqrt{(z-l_{j})(z-r_{j})}}, \qquad \nu_{j} \in I_{j}, \quad z \in \mathbb{C}_{+},
		\]
		where the square roots are normalized as $s_j(z)\sim z-c_j$. Indeed, $\arg g_q=\pi/2$ on the open bands and, by the Herglotz property, changes from $\pi$ to $0$ at most once on each open gap. The crossing point, with an endpoint used if the sign is constant, is $\nu_j$. The factor $i/2$ follows from $g_q(iy)\to i/2$. On $E_N^\circ$, reflectionlessness gives
		\[
		\Im g_{q}=\dfrac{1}{2} \dfrac{1+|m_{-}^{q}|^{2}}{1-|m_{-}^{q}|^{2}} \geq \dfrac{1}{2}.
		\]
		If $\nu_j=e\in\partial I_j$, $g_{q}$ would tend to zero along the $E$, contradicting this bound. Thus
		\[
		\nu_{j} \in I_j^\circ, \qquad \Delta^{-1} = \dfrac{1}{2}-ig_{q}, \qquad \sqrt{\Delta} = c_{e}(z-e)^{1/4}(1+O(|z-e|^{1/2})),
		\]
		with $c_e\ne0$. More precisely, in a slit neighborhood of $e$, the product shows that $(z-e)^{-1/4}\sqrt\Delta$ is analytic in $\sqrt{z-e}$ and nonzero at $e$. From the local factors $W=H_jf_j$, $(z-e)^{1/4}W$ is also analytic in $\sqrt{z-e}$. Thus $d$ and $d^{-1}$ are analytic in this parameter at $e$ and neither vanishes there. The same holds for a local branch of $\log d$, and
		\begin{equation}
			\eta(z)=\eta(e)+O(|z-e|^{1/2}),\qquad \eta(e)\in\mathbb R. \label{3.16}
		\end{equation}
		
		By \eqref{3.14} and $W_{+} = iQ_{\phi_{j,N}+\theta_{j}}W_{-}$,
		\[
		2\Re \eta_{+} = \eta_{+}+\eta_{-} = \phi_{j,N}+\theta_{j}-\phi_{j,N,\theta} \pmod{2\pi}\quad\text{on }I_j^\circ. 
		\]
		Set $\rho(s)=-\partial_{s}\Re \eta_{+}$ and $d\mu(s)=\rho(s)ds$. Then
		\[
		|\mu|(I_{j}) \leq \dfrac{1}{2} \left( \operatorname{Var}_{I_{j}} \phi_{j,N}+ \operatorname{Var}_{I_{j}} \phi_{j,N,\theta}\right) \leq \pi.
		\]
		This gives \eqref{3.15}. The expansion in $\sqrt{z-e}$ also gives
		\[
		u(s) = O(|s-e|^{1/2}), \qquad \rho(s)=O(|s-e|^{-1/2}), \qquad \rho \in L^{3/2}(G_{N}).
		\]
		
		Set $U(x+iy)=\Im\eta(x+i|y|)$ for $y\ne0$ and $U(x)=u(x)$. By the Cauchy--Riemann equations and even reflection,
		\[
		-\partial_{y} U(x,0+)=-\partial_{x}\Re \eta_{+}(x)=\rho(x), \qquad -\Delta U=2\rho(s)ds \quad \text{on } \Omega_{N}.
		\]
		By \eqref{3.16}, $U\in C(\mathbb C)$, $U|_{E_N}=0$ and $U(\infty)=0$. Hence \eqref{3.17} gives $u=(\mathcal K_{\Omega_N}\mu)|_{\mathbb R}$, with the potential extended by zero on $E_N$. Finally, since $\eta_+-\eta_-=2iu$ and $\eta(z)=O(z^{-1})$ at infinity, Cauchy's formula and \eqref{3.16} yield
		\[
		\eta(z) = \dfrac{1}{2\pi i} \int_{G_{N}} \dfrac{\eta_{+}(s)-\eta_{-}(s)}{s-z} ds =  \dfrac{1}{\pi} \int_{G_{N}} \dfrac{u(s)}{s-z}ds.
		\]
		\bigskip
	\end{proof}
	
	\begin{lemma}
		\label{l7} Suppose $S_2<\infty$ and fix $N$. For $q\in\mathcal R(E_N)$, $x\in\mathbb R$ and $z\in\mathbb C_+$, let
		\begin{equation}
			R_{q}(x,z)=i\dfrac{(1-m_{+}^{q}(x,z))(1-m_{-}^{q}(x,z))}{1-m_{+}^{q}(x,z)m_{-}^{q}(x,z)} = i \prod_{j \in \mathcal{J}_{N}} \dfrac{z-\nu_{j}(q;x)}{\sqrt{(z-l_{j})(z-r_{j})}}. \label{3.26}
		\end{equation}
		Here $\nu_j(q;x)$ is the unique zero in $I_j^\circ$ if one exists, and equals $l_j$ or $r_j$ if $R_q$ is positive or negative throughout $I_j^\circ$, respectively. Set
		\[
		Q_{n}(q;x)=Q_{n,N}(q;x)=\sum_{j \in \mathcal{J}_{N}}\left( l_{j}^{n}+r_{j}^{n}-2\nu_{j}(q;x)^{n} \right), \quad n\ge1.
		\]
		Writing $q=u+iv$ for its smooth representative, one has
		\begin{equation}
			Q_{1}=-2u, \qquad Q_{2}=2(v_{x}+v^{2}), \qquad Q_{3}=\dfrac{3}{2}u_{xx}+3vu_{x}-2u^{3}. \label{3.18}
		\end{equation}
		Denote
		\[
		R_{j}=\max(|l_{j}|,|r_{j}|), \qquad M_{n}=2n\sum_{j \in \mathcal{J}_{N}}a_{j}R_{j}^{n-1} \leq 2n(1+S_{0})^{n-1}S_{n-1}.
		\]
		For all $x\in\mathbb R$,
		\begin{equation}
			|q| \leq S_{0}, \qquad |q_{x}| \leq \dfrac{M_{2}}{2}+|q|^{2}, \qquad |q_{xx}|\leq \dfrac{2M_{3}}{3} + 2|q||q_{x}|+\dfrac{4}{3}|q|^{3}. \label{3.20}
		\end{equation}
		In particular,
		\begin{equation}
			\sup_{N,\theta,x,\mathbf t} \left( |q_{N,\theta}(\mathbf t,x)|+|\partial_{x}q_{N,\theta}(\mathbf t,x)| + |\partial_{x}^{2}q_{N,\theta}(\mathbf t,x)| \right) \leq C(S_{0},S_{1},S_{2}). \label{3.21}
		\end{equation}
		Moreover, for $n=1,2,3$,
		\begin{equation}
			\sup_{N\ge L,\,\theta,x,\mathbf t,\gamma}
			\sum_{j\in\mathcal J_N\setminus\mathcal J_L}
			\bigl|l_j^n+r_j^n-2\nu_j(e^{i\gamma}q_{N,\theta}(\mathbf t,\cdot);x)^n\bigr|
			\le 2n\sum_{j\in\mathcal J\setminus\mathcal J_L}a_jR_j^{n-1}
			\xrightarrow[L\to\infty]{}0. \label{3.19}
		\end{equation}
	\end{lemma}
	
	\begin{proof}
		Put $h_\pm(x,z)=i(1+m_\pm^q(x,z))/(1-m_\pm^q(x,z))$. Then
		\[
		\Im h_{\pm}>0, \qquad R_{q}=-\dfrac{2}{h_{+}+h_{-}},
		\]
		so $R_q$ is Herglotz. Reflectionlessness gives the product in \eqref{3.26} by \cite[(4.12)]{FLLZ}. Since $G_N$ is bounded, \cite[Theorem~4.1]{RZ} gives, after rescaling,
		a smooth representative of $q$ with bounded derivatives of all orders. In our normalization the spatial Riccati equations are
		\[
		2zm_{+}^{q}=q+i\partial_{x}m_{+}^{q}+\overline{q}(m_{+}^{q})^{2}, \qquad 2zm_{-}^{q}=\overline{q}-i\partial_{x}m_{-}^{q}+q(m_{-}^{q})^{2}.
		\]
		Using the high-energy Weyl expansion \cite[Section~4]{CG} and comparing coefficients in these identities gives, as $z=iy$ tends to $i\infty$,
		\[
		\begin{array}{c}
			m_{+}^{q}=\dfrac{q}{2z}+\dfrac{iq_{x}}{4z^{2}}+\dfrac{-q_{xx}+|q|^{2}q}{8z^{3}}+O(z^{-4}), \\
			m_{-}^{q}=\dfrac{\overline{q}}{2z}-\dfrac{i\overline{q}_{x}}{4z^{2}}+\dfrac{-\overline{q}_{xx}+|q|^{2}\overline{q}}{8z^{3}}+O(z^{-4}).
		\end{array}
		\]
		Thus,
		\begin{align*}
			\log(R_{q}/i)&=\log(1-m_{+}^{q})+\log(1-m_{-}^{q})-\log(1-m_{-}^{q}m_{+}^{q}) \\
			&=-\dfrac{u}{z}+\dfrac{v_{x}+v^{2}}{2z^{2}}+\dfrac{3u_{xx}+6vu_{x}-4u^{3}}{12z^{3}}+O(z^{-4}).
		\end{align*}
		On the other hand, \eqref{3.26} gives $\log(R_q/i)=\sum_{n=1}^3Q_n/(2nz^n)+O(z^{-4})$. Comparison gives \eqref{3.18}.
		
		For any $\nu \in [l_{j},r_{j}]$, the mean value theorem yields
		\[
		\left| l_{j}^{n}+r_{j}^{n}-2\nu^{n} \right| \leq nR_{j}^{n-1}(|l_{j}-\nu|+|r_{j}-\nu|)=2na_{j}R_{j}^{n-1},
		\]
		so $|Q_n|\le M_n$. Summing over the omitted gaps gives \eqref{3.19}.
		
		Finally, let $q_{\gamma} = e^{i\gamma}q$ with $\gamma \in \mathbb{R}$. Since
		\[
		C_{\gamma}D_{q}C_{\gamma}^{-1}=D_{q_{\gamma}},\qquad  C_{\gamma}=\operatorname{diag}(e^{i\gamma/2},e^{-i\gamma/2}),
		\]
		$q_{\gamma} \in \mathcal{R}(E_{N})$. Thus for $q_{\gamma}=u_{\gamma}+iv_{\gamma}$, $|Q_{n}| \leq M_{n}$ implies
		\begin{align*}
			|u_{\gamma}| &\leq M_{1}/2\le S_{0},\\
			|\partial_{x}v_{\gamma}| & \leq M_{2}/2 + |q|^{2}, \\
			|\partial_{x}^{2}u_{\gamma}| & \leq 2M_{3}/3+2|q||q_{x}|+4|q|^{3}/3.
		\end{align*}
		Taking the supremum over $\gamma$ in each inequality gives \eqref{3.20}. Applying it to the time slices in \eqref{3.22} gives \eqref{3.21}. \bigskip
	\end{proof}
	
	\section{Reconstruction for infinitely many gaps}
	
	\subsection{Uniqueness of the admissible matrix}
	\begin{definition}
		Fix $q_0\in\mathcal R(E)$ and $\theta\in\Theta$. A $2\times2$ matrix function $W$, holomorphic on $\mathbb C\setminus\mathbb R$, is \emph{admissible} for $(q_0,\theta)$ if it satisfies the following conditions.
		\begin{enumerate}[wide=0pt, noitemsep, topsep=2pt]
			\item On each open gap $I_j^\circ$, $W$ has continuous one-sided extensions satisfying
			\[
			W_+(s)=iQ_{\phi_j^{q_0}(s)+\theta_j}W_-(s).
			\]
			\item There exist $q\in\mathcal R(E)$ and a real signed Radon measure $\mu$ on $G$, with $|\mu|(I_j^\circ)\leq\pi$ for all $j$, such that the trace $u=(\mathcal K_\Omega\mu)|_G$ belongs to $L^1(G,\omega(s)\,ds)$ and
			\[
			W(z)=e^{i\eta(z)J}\mathcal{M}_q(z),\qquad
			\eta(z)=\frac{1}{\pi}\int_G\frac{u(s)}{s-z}\,ds,\qquad z\in\mathbb C\setminus\mathbb R.
			\]
		\end{enumerate}
		Here $u$ is extended by zero on $E$.
	\end{definition}
	
	\begin{remark} \label{r2}
		Let $z=x+iy$ with $y>0$ and $(\mathcal K_\Omega|\mu|)(z)<\infty$.
		By the Poisson representation of $g_\Omega(\cdot,t)$ and Tonelli's
		theorem applied to the positive and negative parts of $\mu$,
		\begin{align*}
			\Im\eta(z)
			&=\int_G P_y(x-s)u(s)\,ds\\
			&=\frac1\pi\int_G
			\left(\int_G P_y(x-s)g_\Omega(s,t)\,ds\right)d\mu(t)
			=(\mathcal K_\Omega\mu)(z).
		\end{align*}
	\end{remark}
	
	A nonnegative harmonic function $h$ on $\Omega$ is \emph{quasibounded} if it is the increasing limit of bounded nonnegative harmonic functions on $\Omega$. For a holomorphic function $f$ on $\Omega$, the harmonic majorant criterion \cite[Section~1]{Schiff1979} gives
	\begin{equation}
		f\in\mathcal N_+(\Omega)
		\quad\Longleftrightarrow\quad
		\log^+|f|\leq h\ \text{for some quasibounded harmonic function }h,
		\label{4.1}
	\end{equation}
	where $\log^{+} (t) = \max\{ \log t,0 \}$.
	
	\begin{lemma}\label{l8}
		Suppose that $\Omega$ is a Dirichlet regular Widom domain with DCT. Then for each $(q_0,\theta)$ there is at most one admissible matrix.
	\end{lemma}
	
	\begin{proof}
		Let $W_k=e^{i\eta_kJ}\mathcal{M}_k$ be admissible for $(q_0,\theta)$, where $\mathcal{M}_k=\mathcal{M}_{q_k}$, $k=1,2$. Set
		\[
		V=W_2^{-1}W_1=\mathcal{M}_2^{-1}e^{i(\eta_1-\eta_2)J}\mathcal{M}_1 .
		\]
		Since $\det \mathcal{M}_k=1$ and the two matrices have the same left jump, $V_+=V_-$ on every open gap $I_{j}^{\circ}$. Thus $V$ extends holomorphically to $\Omega$ by Morera's theorem.
		
		For $m_k$ defined by \eqref{3.24} with $q_k$ in place of $q_0$, set $p_k=m_km_k^\#$ and
		\[
		g_k(z)=\frac{i}{2}\frac{1+p_k(z)}{1-p_k(z)}\quad(z\in\mathbb C_+),
		\qquad g_k(\bar z)=\overline{g_k(z)} .
		\]
		By \eqref{3.25}, the boundary values $m_k(\lambda+i0)$ and
		$m_k(\lambda-i0)$ are distinct points of the unit circle for every
		$\lambda\in I_j^\circ$. Thus
		\[
		p_k(\lambda+i0)
		=m_k(\lambda+i0)\overline{m_k(\lambda-i0)}
		\in\mathbb \{ |z|=1 \}\setminus\{1\}.
		\]
		Hence $g_k$ is a Herglotz function with a holomorphic, real-symmetric extension to $\Omega$. By \cite[Theorem~2.1]{BLY}, its lift to the disk has bounded characteristic and no singular inner factor. As $g_k$ has no poles, $g_k\in\mathcal N_+(\Omega)$. Since $\det \mathcal{M}_k=1$ and $(1-p_k)^{-1}=1/2-ig_k$ in $\mathbb C_+$, \eqref{4.1} and symmetry give a quasibounded harmonic function $H_k$ such that
		\begin{equation}
			\log\norm{\mathcal{M}_k^{-1}}=\log\norm{\mathcal{M}_k}
			\leq \tfrac12\log(2+4|g_k|)\leq H_k
			\qquad(z\in\mathbb C\setminus\mathbb R).
			\label{4.2}
		\end{equation}
		
		Write $H=H_1+H_2$, $\nu=|\mu_1|+|\mu_2|$ and $P=\mathcal K_\Omega\nu$. Since $\nu$ is a positive, $P$ is
		superharmonic in $\Omega$. The Widom condition gives
		\[
		P(i)\leq 2\sum_j\sup_{t\in I_j^\circ}g_\Omega(i,t)<\infty .
		\]
		Remark \ref{r2} and \eqref{4.2} give
		\[
		|\Im(\eta_1-\eta_2)|\leq P,\qquad \log^+\norm V\leq H+P
		\qquad(z\notin\mathbb R).
		\]
		The nonnegative subharmonic function $v=(\log^+\norm V-H)_+$ satisfies $v\leq P$ throughout $\Omega$, by the mean inequalities for $v$ and $P$ on circles centered on $G$. Exhaust $\Omega$ by increasing smooth domains $D_n\Subset\Omega$ containing $i$, and let $h_n$ be harmonic on $D_n$ with boundary values $v$. Comparison gives
		\[
		0\leq v\leq h_n\leq P,\qquad h_n\leq h_{n+1}\quad\text{on }D_n .
		\]
		By Harnack's convergence theorem the limit is a nonnegative harmonic minorant of the Green potential $P$. It vanishes by \eqref{2.1}. Thus $v=0$, and each entry of $V$ belongs to $\mathcal N_+(\Omega)$ by \eqref{4.1}.
		
		By boundary uniqueness, $1-p_k\ne0$ almost everywhere on $E$, so the
		reflectionless identity gives $\mathcal M_{k,+}=\mathcal M_{k,-}$ there.
		The Cauchy boundary formula gives $\eta_{k,+}-\eta_{k,-}=2iu_k=0$
		almost everywhere on $E$. In particular $V_+=V_-$ there and, by \eqref{4.2},
		\[
		\norm{V_\pm}\leq \sqrt{(2+4|g_{1,+}|)(2+4|g_{2,+}|)}
		\leq 2+2|g_{1,+}|+2|g_{2,+}|.
		\]
		Here $g_{k,+}\in i[0,\infty)$ almost everywhere on $E$, so the Herglotz representation gives
		\[
		\int_E\frac{|g_{k,+}(s)|}{1+s^2}\,ds
		\leq \pi\Im g_k(i).
		\]
		For $\xi_*\in G$, $|s-\xi_*|^{-2}\leq C_{\xi_*}(1+s^2)^{-1}$ on $E$. We may therefore apply DCT for $V$ entrywise:
		\[
		V(z)-V(\xi_*)=\frac{1}{2\pi i}\int_E(V_+(s)-V_-(s))
		\left(\frac{1}{s-z}-\frac{1}{s-\xi_*}\right)\,ds=0 .
		\]
		Finally, $|s-iy|^{-1}\leq\omega(s)$ for $y\geq1$. Dominated convergence gives $\eta_k(iy)\to0$, and the Weyl normalization gives $\mathcal{M}_k(iy)\to I_2$. Hence $V\equiv I_2$, namely $W_{1} = W_{2}$. \bigskip
	\end{proof}
	
	\subsection{Existence of the admissible matrix}
	Set $X=L^1(G,\omega(s)\,ds)$ and
	\[
	\mathcal U=\bigcup_N\mathcal U_N
	=\bigcup_N
	\left\{
	\mathcal K_{\Omega_N}\mu :
	\begin{array}{l}
		\mu\text{ is a real signed measure on }G_N,\\
		|\mu|(I_j^\circ)\leq\pi,\quad j\in\mathcal J_N
	\end{array}
	\right\},
	\]
	where $\mathcal K_{\Omega_N}\mu$ extended by zero on $G\setminus G_N$.
	
	\begin{lemma}\label{l9}
		Suppose that $\Omega$ is Dirichlet regular, $S_0<\infty$, and
		$\mathcal P_1(E)<\infty$. Then $\mathcal U$ is relatively compact
		in $X$. If $N_k\to\infty$ and
		$u_k=\mathcal K_{\Omega_{N_k}}\mu_k\in\mathcal U_{N_k}$, then,
		after passage to a subsequence,
		\[
		u_k\longrightarrow u= \mathcal K_\Omega\mu \quad\text{in }X
		\]
		for some real signed Radon measure $\mu$ on $G$ satisfying
		$|\mu|(I_j^\circ)\leq\pi$ for every $j$. Moreover,
		\[
		\lim_{y\to\infty}\sup_{u\in\mathcal U}
		\left|\frac1\pi\int_G\frac{u(s)}{s-iy}\,ds\right|=0.
		\]
	\end{lemma}
	
	\begin{proof}
		Let $v_{\omega,N}=\mathcal K_{\Omega_N}(\omega|_{G_N})$.
		By symmetry of the Green kernel, \eqref{2.15}, and Tonelli's theorem,
		\begin{equation}\label{4.4}
			\begin{aligned}
				\left\|\frac1\pi\int_{G_N\setminus G_L}
				g_{\Omega_N}(\cdot,t)\,d\mu(t)\right\|_X
				&\leq \int_{G_N\setminus G_L}v_{\omega,N}(t)\,d|\mu|(t)\\
				&\leq \pi\sum_{j\notin\mathcal J_L}\sup_{I_j}v_\omega
				\xrightarrow[L\to\infty]{}0
			\end{aligned}
		\end{equation}
		uniformly in $N$ and $\mu$. The same estimate holds for $g_\Omega$.
		
		Fix $j$ and extend the Green kernels by zero at the endpoints of $I_j$.
		By Lemma~\ref{l1}, $v_{\omega,N},v_\omega\in C(I_j)$ and
		$v_{\omega,N}\uparrow v_\omega$. Thus Dini's theorem gives
		\begin{equation}\label{4.3}
			\sup_{t\in I_j}
			\|g_\Omega(\cdot,t)-g_{\Omega_N}(\cdot,t)\|_X
			=\pi\sup_{t\in I_j}\bigl(v_\omega(t)-v_{\omega,N}(t)\bigr)
			\longrightarrow0.
		\end{equation}
		
		Suppose $N_k\to\infty$. Let $\mu_k^+$ and $\mu_k^-$ denote the positive and negative parts of $\mu_k$. Regard
		$\mu_{k,j}^\pm:=\mu_k^\pm|_{I_j^\circ}$ as measures on $I_j$
		with zero endpoint masses. Passing to a diagonal subsequence, we have
		\[
		\mu_{k,j}^\pm\stackrel{*}{\rightharpoonup}\nu_j^\pm
		\quad\text{on }I_j,\qquad
		\nu_j^+(I_j)+\nu_j^-(I_j)
		=\lim_{k\to\infty}|\mu_k|(I_j^\circ)\leq\pi .
		\]
		Define $\mu$ on $G$ by
		$\mu|_{I_j^\circ}=(\nu_j^+-\nu_j^-)|_{I_j^\circ}$.
		Then $|\mu|(I_j^\circ)\leq\pi$.
		
		Fix either sign and write $\sigma_k=\mu_{k,j}^\pm$,
		$\sigma=\nu_j^\pm$. Put
		\[
		w_k(s)=\frac1\pi\int_{I_j}g_\Omega(s,t)\,d\sigma_k(t),
		\qquad
		w(s)=\frac1\pi\int_{I_j}g_\Omega(s,t)\,d\sigma(t).
		\]
		Since $v_\omega\in C(I_j)$,
		\[
		\|w_k\|_X=\int_{I_j}v_\omega\,d\sigma_k
		\longrightarrow \int_{I_j}v_\omega\,d\sigma=\|w\|_X .
		\]
		For each $s\in G$ and $M>0$,
		$\min\{g_\Omega(s,\cdot),M\}\in C(I_j)$. Hence
		\[
		\liminf_{k\to\infty}w_k(s)
		\geq \frac1\pi\int_{I_j}\min\{g_\Omega(s,t),M\}\,d\sigma(t)
		\xrightarrow[M\to\infty]{}w(s).
		\]
		Thus $(w-w_k)_+\to0$ a.e. and $0\leq(w-w_k)_+\leq w\in X$.
		By dominated convergence,
		\[
		\|w_k-w\|_X
		=\|w_k\|_X-\|w\|_X+2\|(w-w_k)_+\|_X
		\longrightarrow0 .
		\]
		
		Since $g_\Omega(\cdot,l_j)=g_\Omega(\cdot,r_j)=0$, subtraction of the
		two limits and \eqref{4.3} give, for each fixed $j$,
		\[
		\frac1\pi\int_{I_j^\circ} g_{\Omega_{N_k}}(\cdot,t)\,d\mu_k(t)
		\xrightarrow[k\to\infty]{X}
		\frac1\pi\int_{I_j} g_\Omega(\cdot,t)\,d(\nu_j^+-\nu_j^-)(t)
		=
		\frac1\pi\int_{I_j^\circ} g_\Omega(\cdot,t)\,d\mu(t).
		\]
		Summing over $j\in\mathcal J_L$ and applying \eqref{4.4} to both tails,
		we obtain
		\[
		\limsup_{k\to\infty}\|u_k-\mathcal K_\Omega\mu\|_X
		\leq 2\pi\sum_{j\notin\mathcal J_L}\sup_{I_j}v_\omega
		\xrightarrow[L\to\infty]{}0 .
		\]
		If the indices $N_k$ are bounded, pass to a constant subsequence
		$N_k=N$ and apply the same argument on the finitely many gaps of
		$G_N$, with $g_{\Omega_N}$ in place of $g_\Omega$.
		This proves relative compactness.
		
		Finally, let $T_yu=\pi^{-1}\int_G u(s)/(s-iy)\,ds$. Since
		$|s-iy|^{-1}\leq\omega(s)$ for $y\geq1$, dominated convergence gives
		\[
		\sup_{y\geq1}\|T_y\|_{X^*}\leq\frac1\pi,\qquad
		T_yu\longrightarrow0\quad(u\in X).
		\]
		The convergence is uniform on the compact closure of $\mathcal U$. \bigskip
	\end{proof}
	
	\begin{lemma}\label{l10}
		Let $\Omega$ be a Dirichlet regular Widom domain with DCT, and suppose that $S_2<\infty$ and $\mathcal P_1(E)<\infty$. Fix $q_0\in\mathcal R(E)$, write $W_{N,\theta}(z)=W_N(\theta,z)$ and put $q_{N,\theta}(x)=F_N(\theta+x\alpha_N)$. For each $\theta\in\Theta$ there exists a unique admissible matrix $W_\theta$ for $(q_0,\theta)$, and
		\[
		W_{N_n,\theta^{(n)}}\longrightarrow W_\theta
		\qquad\text{whenever }N_n\to\infty,\quad\theta^{(n)}\to\theta,
		\]
		locally uniformly on $\mathbb C\setminus\mathbb R$. The map $\theta\mapsto W_\theta$ is continuous in this topology. Moreover,
		\[
		F(\theta):=\lim_{y\to\infty}2iy(W_\theta(iy))_{21},\qquad
		\norm{F_N-F}_{C(\Theta)}\longrightarrow0,\qquad W_0=\mathcal{M}_{q_0}.
		\]
	\end{lemma}
	
	\begin{proof}
		Fix $N_n\to\infty$ and $\theta^{(n)}\to\theta$, and abbreviate the pair of indices by $n$. By \eqref{3.22} and Lemmas~\ref{l6} and~\ref{l7},
		\[
		q_n\in\mathcal R(E_{N_n}),\qquad
		\sup_n(\norm{q_n}_\infty+\norm{q_n'}_\infty+\norm{q_n''}_\infty)<\infty,
		\qquad W_n=e^{i\eta_nJ}\mathcal{M}_{q_n},
		\]
		and $u_n=\mathcal K_{\Omega_{N_n}}\mu_n$ on $G$ and $\eta_n$ is its Cauchy transform. By Arzel\`a--Ascoli and Lemma~\ref{l9}, after passage to a subsequence,
		\[
		q_n\to q\ \text{in } C^1_{\mathrm{loc}}(\mathbb R),\qquad
		u_n\to u=\mathcal K_\Omega\mu\ \text{in }X,\qquad |\mu|(I_j^\circ)\leq\pi .
		\]
		Then the half-line Weyl functions converge locally uniformly off the real line \cite[Proposition~2.3 and Section~3.4]{EKT}. Consequently,
		\begin{equation}
			\mathcal{M}_{q_n}\to \mathcal{M}_q,\qquad
			\eta_n\to\eta=\frac1\pi\int_G\frac{u(s)}{s-z}\,ds,\qquad
			W_n\to W=e^{i\eta J}\mathcal{M}_q
			\label{4.5}
		\end{equation}
		locally uniformly on $\mathbb C\setminus\mathbb R$.
		
		On $\mathbb C_+$ write $a_n=m_-^{q_n}$ and $b_n(z)=\overline{m_+^{q_n}(\bar z)}$, with respective limits $a,b$. Thus
		\[
		\mathcal{M}_q=(1-ab)^{-1/2}\begin{pmatrix}1&a\\b&1\end{pmatrix},
		\]
		where the square root has positive real part. Local uniform convergence of bounded Schur functions implies weak$^*$ convergence of their boundary values in $L^\infty(\mathbb R)$, by the Poisson representation. Since $a_n=\overline{b_n}$ almost everywhere on $E\subset E_{N_n}$, the same holds for $a,b$. For $\lambda\in I_j^\circ$ and $d=\operatorname{dist}(\lambda,\partial I_j)$, the spectral theorem gives, for all sufficiently large $n$,
		\[
		\norm{(D_{q_n}-\lambda)f}_2\geq d\norm f_2
		\qquad \bigl(f\in H^1(\mathbb R;\mathbb C^2)\bigr).
		\]
		Passing to the limit by dominated convergence gives the same bound for
		$D_q$, hence $\sigma(D_q)\subset E$. Conversely, let $g_q=\frac{i}{2}(1+ab)/(1-ab)$. As in the finite-gap
		case, the reflectionless identity gives
		\[
		\Im g_q(s+i0)\geq \dfrac12
		\qquad\text{for a.e. }s\in E.
		\]
		By Dirichlet regularity and \cite[Section~2.1]{BLY},
		$\operatorname{supp}(\mathbf 1_E\,ds)=E$. Thus the preceding lower bound
		and Stieltjes inversion give
		\[
		E=\operatorname{supp}(\mathbf 1_E\,ds)
		\subset \operatorname{supp}\mu_{g_q}=\sigma(D_q),
		\]
		where $\mu_{g_q}$ is the representing measure of $g_q$.
		Hence $q\in\mathcal R(E)$.
		
		Fix a compact interval $A\subset I_j^\circ$. The full-line Green matrix at $x=0$ is
		\[
		\mathcal G_n(z)=\frac{i}{2(1-a_nb_n)}
		\begin{pmatrix}1+a_nb_n&2b_n\\2a_n&1+a_nb_n\end{pmatrix},
		\qquad h_n=\dfrac12-i(\mathcal G_n)_{11}.
		\]
		Its matrix-valued Herglotz measure is supported on $E_{N_n}$, and $\sup_n\norm{\mathcal G_n(i)}<\infty$. On a fixed neighborhood of $A$ the Herglotz kernel $(s-z)^{-1}-s/(1+s^2)$ is $O((1+s^2)^{-1})$ uniformly for $s\in E_{N_n}$ and all sufficiently large $n$. The Herglotz representation gives local boundedness there, and the identity theorem identifies the unique limit with $\mathcal G_q$. Since $(\mathcal G_n)_{11}$ is real on $A$, one has $|h_n|\geq1/2$ there. On a smaller neighborhood the analytic continuations
		\[
		a_n=\frac{(\mathcal G_n)_{21}}{ih_n},\qquad
		b_n=\frac{(\mathcal G_n)_{12}}{ih_n}
		\]
		therefore converge locally uniformly, with all derivatives. As $h_n=(1-a_nb_n)^{-1}$,
		\[
		(\mathcal{M}_{q_{n}})_{\pm} \to (\mathcal{M}_{q})_{\pm} \quad \text{on } A.
		\]
		
		On $I_j^\circ$ write $d\mu_n(s)=\rho_n(s)\,ds$. By \eqref{3.15} and the choice of phase in \eqref{3.2},
		\[
		\rho_n=-\dfrac12\bigl(\phi_{j,N_n}'-\phi_{j,N_n,\theta^{(n)}}'\bigr),
		\qquad
		\phi_{j,N_n,\theta^{(n)}}'
		=\dfrac1{2i}\left(\frac{a_n'}{a_n}-\frac{b_n'}{b_n}\right).
		\]
		Since $|a_n|=|b_n|=1$ on the gap, \eqref{3.28} and the convergence just proved imply local uniform convergence of $\rho_n$ on $I_j^\circ$. Choose $\chi\in C_c^\infty(I_j^\circ)$ equal to one on a neighborhood of $A$, and set
		\[
		L_n(z)=-\frac i\pi\int_{I_j}\chi(s)\rho_n(s)\log(z-s)\,ds
		\quad (z\in\mathbb C_+),\qquad 0<\arg(z-s)<\pi .
		\]
		The local integrability of the logarithm gives uniform convergence of $L_n$ on compact subsets up to the upper side of $A$. Put $U_n(x,y)=\Im\eta_n(x+i|y|)$ and extend $\Im L_n$ evenly across $A$. On any sufficiently small disk $B$ centered on $A$ with $\chi=1$ on $B\cap\mathbb R$,
		\[
		-\Delta U_n=-\Delta\Im L_n=2\rho_n(s)\,ds,\qquad
		\norm{U_n-U_m}_{L^1(B)}\leq C_B\norm{u_n-u_m}_X .
		\]
		The integral of the Poisson kernel over $B$ is bounded for $s$ in a bounded interval and is $O(|s|^{-2})$ at infinity, which gives the inequality. Hence $h_{n}=U_n-\Im L_n$ is a Cauchy sequence of harmonic functions in $L^1(B)$ and
		for every disk $B'\Subset B$ we have
		\[
		\sup_{B'}|h_n-h_m|+\sup_{B'}|\nabla h_n-\nabla h_m|
		\leq C_{B',B}\int_B |h_n-h_m|\,dA
		\xrightarrow[n,m\to\infty]{}0 .
		\]
		By the Cauchy--Riemann equations and \eqref{4.5}, for every compact
		interval $A\subset I_j^\circ$,
		\[
		\lim_{n,m\to\infty}\sup_{s\in A}
		\left|(\eta_n-L_n)_+(s)-(\eta_m-L_m)_+(s)\right|=0.
		\]
		Adding $L_n$ and using symmetry on the lower side gives
		\[
		(\eta_n)_\pm\to\eta_\pm,\qquad (W_n)_\pm\to W_\pm
		\quad\text{locally uniformly on }I_j^\circ .
		\]
		Passing to the limit in \eqref{3.6} and using \eqref{3.28}, we obtain
		$W_+=iQ_{\phi_j^{q_0}+\theta_j}W_-$ on $I_j^\circ$. This proves admissibility of $W$.
		
		By Lemma~\ref{l8}, $W=W_\theta$, so every subsequence has the same
		limit point and then $W_n\to W_\theta$ locally uniformly on
		$\mathbb C\setminus\mathbb R$. If $\theta^{(n)}\to\theta$, choose $N_n\geq n$ so that
		\[
		\left|W_{\theta^{(n)}} - W_{N_{n},\theta^{(n)}}\right| \leq \dfrac{1}{n}
		\]
		on the $n$th compact set of an increasing exhaustion of $\mathbb C\setminus\mathbb R$. The same convergence gives
		\[
		|W_{\theta^{(n)}} - W_{\theta}| \leq |W_{\theta^{(n)}} - W_{N_{n},\theta^{(n)}}|+|W_{N_{n},\theta^{(n)}}-W_{\theta}| \to 0.
		\]
		Namely, $W_{\theta^{(n)}}\to W_\theta$ locally uniformly.
		
		Finally, for $y>0$, the Riccati equation for $m_-^q$ gives
		\[
		m_-^q(x,iy)
		=-i\int_{-\infty}^{x}e^{-2y(x-s)}
		\left(\overline{q(s)}+q(s)m_-^q(s,iy)^2\right)\,ds .
		\]
		Since $|m_-^q|\leq1$, we have $|m_-^q(x,iy)|\leq \norm{q}_{\infty}/y$. Hence, by continuity of $q$ at zero,
		\[
		2iy\,m_-^q(0,iy)
		=\int_{-\infty}^{0}2ye^{2ys}\overline{q(s)}\,ds+O(y^{-2})
		\longrightarrow\overline{q(0)} .
		\]
		Similarly, for the right Weyl coordinate in $\mathbb C_+$,
		\[
		2iy\,m_+^q(0,iy)
		=\int_0^\infty 2ye^{-2ys}q(s)\,ds+O(y^{-2})
		\longrightarrow q(0).
		\]
		Since $(\mathcal{M}_q)_{21}=m_+^q(\mathcal{M}_q)_{22}$ in $\mathbb C_+$,
		$\mathcal{M}_q(iy)\to I_2$, and $\eta(iy)\to0$, we obtain
		\[
		F(\theta):=\lim_{y\to\infty}2iy\bigl(W_\theta(iy)\bigr)_{21}
		=\lim_{y\to\infty}2iy\bigl(\mathcal{M}_q(iy)\bigr)_{21}
		=q(0).
		\]
		
		Let $N_k\to\infty$ and $\theta^{(k)}\to\theta$. From every subsequence
		we can extract a further subsequence along which
		\[
		q_{N_k,\theta^{(k)}}\longrightarrow q
		\quad\text{in }C^1_{\mathrm{loc}}(\mathbb R).
		\]
		For each such limit, the preceding identity gives $q(0)=F(\theta)$.
		It follows that along the original sequence,
		\begin{equation}\label{4.6}
			F_{N_k}(\theta^{(k)})
			=q_{N_k,\theta^{(k)}}(0)
			\longrightarrow F(\theta).
		\end{equation}
		
		Set
		\[
		d_L=\sup_{N,M\geq L}\sup_{\theta\in\Theta}
		|F_N(\theta)-F_M(\theta)|,\qquad L\geq1.
		\]
		For each $L$, choose $N_L,M_L\geq L$ and $\theta_L\in\Theta$ such that
		\[
		0\leq d_L
		\leq |F_{N_L}(\theta_L)-F_{M_L}(\theta_L)|+\frac1L .
		\]
		By compactness of $\Theta$, there is a subsequence
		$\theta_{L_j}\to\theta$. Applying \eqref{4.6} twice gives
		\[
		F_{N_{L_j}}(\theta_{L_j})\longrightarrow F(\theta),
		\qquad
		F_{M_{L_j}}(\theta_{L_j})\longrightarrow F(\theta).
		\]
		Thus $d_{L_j}\to0$, and monotonicity gives $d_L\downarrow0$.
		For every $M\geq N$,
		\[
		\sup_{\theta\in\Theta}|F_N(\theta)-F_M(\theta)|\leq d_N .
		\]
		Letting $M\to\infty$ at each $\theta$ and then taking the supremum,
		we obtain
		\[
		\|F_N-F\|_{C(\Theta)}\leq d_N\longrightarrow0 .
		\]
		In particular, $F\in C(\Theta)$.
		
		Finally, taking $q=q_0$ and $\mu=0$ in the definition of
		admissibility, we see from \eqref{3.2} that $\mathcal{M}_{q_0}$ is admissible
		for $(q_0,0)$. By Lemma~\ref{l8}, $W_0=\mathcal{M}_{q_0}$.
	\end{proof}
	
	\subsection{Proof of Theorem \ref{t1}}
	
	Suppose that the hypotheses of Theorem~\ref{t1} hold. We use the reconstruction function $F$ of Lemma~\ref{l10}.
	
	\begin{lemma}\label{l11}
		Let $q_0\in\mathcal R(E)$ and
		\[
		q(\mathbf t,x)=F(\vartheta),\qquad
		\vartheta=x\alpha+\sum_{r=2}^m t_r\beta^{(r)} .
		\]
		There exist $F_0,\ldots,F_m,F_{t_2},\ldots,F_{t_m}\in C(\Theta)$, with $F_0=F$, such that
		\begin{equation}\label{4.10}
			\partial_x^j q(\mathbf t,x)=F_j(\vartheta)\quad(0\le j\le m),\qquad
			\partial_{t_r}q(\mathbf t,x)=F_{t_r}(\vartheta)\quad(2\le r\le m).
		\end{equation}
		Each $x\mapsto F(\theta+x\alpha)$ belongs to $C^m(\mathbb R)$, and
		\begin{equation}\label{4.13}
			\partial_{t_r}q=\mathcal X_r(q)\quad(2\le r\le m),\qquad
			q(\mathbf0,x)=q_0(x)\quad\text{for a.e. }x\in\mathbb R.
		\end{equation}
		Moreover, for $k\in\mathbb N_0$ and $\boldsymbol\ell=(\ell_2,\ldots,\ell_m)\in\mathbb N_0^{m-1}$ with $k+\sum_{r=2}^m r\ell_r\le m$, there exists $F_{k,\boldsymbol\ell}\in C(\Theta)$ such that
		\[
		\partial_x^k\partial_{\mathbf t}^{\boldsymbol\ell}q(\mathbf t,x)=F_{k,\boldsymbol\ell}(\vartheta).
		\]
	\end{lemma}
	
	\begin{proof}
		Put $F_{r,N}(\theta)=\partial_x^r q_{N,\theta}(\mathbf0,0)$ for $0\le r\le m$, so that $F_{0,N}=F_N$. Fix $\theta\in\Theta$ and let $N_n\to\infty$, $\theta_n\to\theta$. For $p_n(x)=q_{N_n,\theta_n}(\mathbf0,x)$ and $p_\theta(x)=F(\theta+x\alpha)$, by Lemma~\ref{l10} and coordinatewise convergence of the frequencies,
		\[
		p_n(x)=F_{N_n}(\theta_n+x\alpha_{N_n})\longrightarrow p_\theta(x)
		\]
		locally uniformly in $x$. By \eqref{3.21} the sequence is relatively compact in $C^1_{\rm loc}(\mathbb R)$. Every limit has value $p_\theta(x)$ at each $x$, so
		\begin{equation}\label{4.7}
			p_n\longrightarrow p_\theta\quad\text{in } C^1_{\rm loc}(\mathbb R).
		\end{equation}
		In particular, $F_{1,N_n}(\theta_n)\to F_1(\theta):=p_\theta'(0)$.
		
		Write $p_n=u_n+iv_n$. At $x=0$, the third trace identity in \eqref{3.18} reads
		\begin{equation}\label{4.8}
			Q_{3,N_n}(p_n)=\sum_{j\in\mathcal J_{N_n}}
			(l_j^3+r_j^3-2\nu_j(p_n;0)^3)
			=\tfrac32u_n''+3v_nu_n'-2u_n^3 .
		\end{equation}
		The Weyl functions, and hence $R_{p_n}(0,\cdot)$, converge locally uniformly in $\mathbb C_+$ by \eqref{4.7}. The Herglotz measures of $R_{p_n}$ are supported on $E_{N_n}$. The Herglotz representation and boundedness of $R_{p_n}(0,i)$ therefore give local uniform convergence across that compact subinterval of $I_{j}^{\circ}$. Since $R_p(0,\cdot)$ is strictly increasing in each gap $I_{j}^{\circ}$,
		\[
		\nu_j(p;0)=l_j+\int_{l_j}^{r_j}\boldsymbol{1}_{\{R_p(0,s)<0\}}\,ds.
		\]
		Consequently $\nu_j(p_n;0)\to\nu_j(p_\theta;0)$ for each $j$, including when the limit is an endpoint. The same argument applies to $ip_n$. For $1\le d\le m+1$ and $\nu\in I_j$,
		\[
		|l_j^d+r_j^d-2\nu^d|\le 2d a_jR_j^{d-1},\qquad
		\sum_j a_jR_j^{d-1}\le(1+S_0)^{d-1}S_{d-1}<\infty .
		\]
		Thus $Q_{d,N_n}(p_n)$ and $Q_{d,N_n}(ip_n)$ converge to limits determined by $p_\theta$. Applying \eqref{4.8} to $p_n$ and $ip_n$ and using \eqref{4.7}, we obtain
		\[
		F_{2,N_n}(\theta_n)=p_n''(0)\longrightarrow F_2(\theta).
		\]
		
		For $2<r\le m$, the same argument applies to $Q_{d,N_{n}}$ of order $d=r+1$. Indeed, in the two Riccati expansions used in Lemma~\ref{l7}, the coefficients of $z^{-d}$ have leading derivative terms $2^{-d}i^{d-1}p^{(d-1)}$ and $2^{-d}(-i)^{d-1}\overline{p}^{(d-1)}$, respectively. This follows by induction since only the linear derivative term raises the order. Comparing the coefficient of $z^{-d}$ in $\log(R_p/i)$ with the finite product \eqref{3.26} gives
		\[
		Q_{d,N}(p)=-d\,2^{2-d}\Re\!\left(i^{d-1}p^{(d-1)}\right)
		+P_d\bigl(p,\bar p,\ldots,p^{(d-2)},\bar p^{(d-2)}\bigr),
		\]
		where $P_d$ is a universal polynomial. Applying this identity to $p_n$ and $ip_n$ recovers both components of $p_n^{(r)}(0)$. Induction now gives
		\[
		F_{r,N_n}(\theta_n)\longrightarrow F_r(\theta),\qquad 0\le r\le m,\quad F_0=F .
		\]
		This holds for every sequence as above. The compactness argument of Lemma~\ref{l10} gives $F_r\in C(\Theta)$ and $\|F_{r,N}-F_r\|_{C(\Theta)}\to0$. For $2\le r\le m$, Lemma~\ref{l5} gives
		\begin{equation}\label{4.11}
			F_{t_r,N}(\theta):=\partial_{t_r}q_{N,\theta}(\mathbf0,0)
			=\left.\mathcal X_r\bigl(q_{N,\theta}(\mathbf0,\cdot)\bigr)\right|_{x=0}.
		\end{equation}
		The right-hand side is a polynomial in $F_{j,N},\overline{F_{j,N}}$ for $0\le j\le r$. Substituting $F_j,\overline{F_j}$ gives $F_{t_r}\in C(\Theta)$ and $F_{t_r,N}\to F_{t_r}$ uniformly on $\Theta$.
		
		For bounded $s$, $s\beta_N^{(r)}\to s\beta^{(r)}$ uniformly in the product metric by Lemma~\ref{l3}, $1\le r\le m$. Integrating the finite-gap phase identities and passing to the limit gives
		\[
		\begin{aligned}
			F(\theta+h\beta^{(r)})-F(\theta)
			&=\int_0^h F_{t_r}(\theta+s\beta^{(r)})\,ds &&(2\le r\le m),\\
			F_j(\theta+h\alpha)-F_j(\theta)
			&=\int_0^h F_{j+1}(\theta+s\alpha)\,ds &&(0\le j<m).
		\end{aligned}
		\]
		Differentiating at $h=0$ gives the directional derivatives at every
		$\theta\in\Theta$. Taking $\theta=\vartheta$ yields \eqref{4.10}. Since $\mathcal X_r$ is
		a differential polynomial of order $r$, passing to the limit in
		\eqref{4.11} and using \eqref{4.10} yields
		\[
		\partial_{t_r}q(\mathbf t,x)
		=F_{t_r}(\vartheta)
		=\mathcal X_r\bigl(q(\mathbf t,\cdot)\bigr)(x),
		\qquad 2\le r\le m.
		\]
		
		Write $q_j=\partial_x^j q=F_j(\vartheta)$ for $0\le j\le m$. For $j+r\le m$, integrate the $r$th equation in $t_r$ and differentiate $j$ times in $x$ under the integral. The integrands are continuous, hence
		\[
		\partial_{t_r}q_j=\partial_x^j\mathcal X_r(q),\qquad
		\partial_{t_r}\overline{q_j}=\overline{\partial_x^j\mathcal X_r(q)}.
		\]
		Each $t_r$-derivative thus raises the largest spatial index in a differential polynomial by at most $r$. Induction on $\sum_r\ell_r$ gives $\partial_x^k\partial_{\mathbf t}^{\boldsymbol\ell}q$ as a polynomial in $q_j,\overline{q_j}$ with $j\le k+\sum_r r\ell_r\le m$. Substituting $F_j,\overline{F_j}$ gives $F_{k,\boldsymbol\ell}\in C(\Theta)$.
		
		Finally, $\mathcal{M}_{q_0}$ is admissible at $\theta=0$ with $\mu=0$. By uniqueness and the limiting factorization \eqref{4.5},
		\[
		\mathcal{M}_{q_0}=W_0=e^{i\eta J}\mathcal{M}_{p_0},\qquad p_0(x)=F(x\alpha).
		\]
		Thus the two half-line Weyl functions of $p_0$ agree with those of $q_0$. Inverse uniqueness for the Dirac operator gives $p_0=q_0$ a.e.\ on $\mathbb R$ \cite[Section~5]{CG}. \bigskip
	\end{proof}
	
	\begin{lemma}\label{l12}
		Let $r\ge2$ and let $u,v$ be classical solutions of $q_t=\mathcal X_r(q)$ on
		$[-T,T]\times\mathbb R$. Suppose that $\partial_x^k u,\partial_x^k v$ are continuous
		and bounded there for $0\le k\le r$. If $u(0,\cdot)=v(0,\cdot)$, then $u=v$.
	\end{lemma}
	
	\begin{proof}
		For $q=u,v$, let $M_q$ be the Weyl matrix
		$\left(\begin{smallmatrix}1&m\\m^\#&1\end{smallmatrix}\right)$
		defined by \eqref{3.24} at $(t,x)$, and write
		$M_0=M_u(0,x,z)=M_v(0,x,z)$. Let
		\[
		\partial_tT_q=T_qV_r(q),\qquad T_q(0,x,z)=I_2.
		\]
		By zero curvature,
		\[
		\partial_xT_q=T_qU_q(t)-U_q(0)T_q.
		\]
		For fixed nonreal $z$, $V_r(q)$ is bounded on the strip, so the
		time equation gives a bound on $T_q$ uniform in $x$; hence right
		multiplication by $T_q$ preserves square integrability at either end.
		Then for $z\in\mathbb C_+$, the first and second rows of $M_0T_q$ are
		proportional to the values at $x$ of the Weyl row solutions at time $t$
		on $\mathbb R_-$ and $\mathbb R_+$, respectively.
		
		Fix a nonreal ray. The spatial Riccati expansion gives, uniformly on
		the strip,
		\[
		M_q=I_2+O(|z|^{-1}),\qquad
		M_q^{-1}JM_q=\sum_{k=0}^r R_k(q)z^{-k}+O(|z|^{-r-1}).
		\]
		Indeed, truncation at order $r$ leaves a residual $O(|z|^{-r})$
		in the Riccati equations, and
		$\|(\partial_x\pm2iz)^{-1}\|_{\infty\to\infty}=(2|\Im z|)^{-1}$.
		The coefficients are determined by $R_x=[R,U_q]$, $R^2=I_2$ and
		$R_0=J$. Hence \eqref{3.12} gives
		\[
		M_qV_r(q)=-iz^rJM_q+O(|z|^{-1}).
		\]
		
		Put $Y_q=e^{itz^rJ}M_0T_q$. By uniqueness of the half-line Weyl
		solutions, we have
		\[
		Y_q=
		\begin{pmatrix}
			(Y_q)_{11}&0\\
			0&(Y_q)_{22}
		\end{pmatrix}M_q,
		\]
		and
		\[
		\partial_t(Y_q)_{jj}
		=\bigl(iz^rJ_{jj}+(M_qV_r(q))_{jj}\bigr)(Y_q)_{jj}
		=O(|z|^{-1})(Y_q)_{jj}.
		\]
		As $(Y_q)_{jj}=1$ at $t=0$, integration gives
		$(Y_q)_{jj}=1+O(|z|^{-1})$. The row relation now yields
		\[
		Y_q=I_2+O(|z|^{-1}),\qquad
		T_u^{-1}T_v=Y_u^{-1}Y_v=I_2+O(|z|^{-1}).
		\]
		
		For fixed $t,x$, $T_u^{-1}T_v$ is entire in $z$ and has growth at
		most $\exp(C_T(1+|z|)^r)$ by the time equation. Applying the
		Phragm\'en--Lindel\"of principle between successive nonreal rays
		of angular separation less than $\pi/r$, and then Liouville's
		theorem, gives $T_u^{-1}T_v=I_2$. Write $T=T_u=T_v$. Subtracting the two spatial equations and using
		$U_u(0)=U_v(0)$, we obtain
		\[
		0=\partial_xT_u-\partial_xT_v
		=T\bigl(U_u(t)-U_v(t)\bigr)
		-\bigl(U_u(0)-U_v(0)\bigr)T
		=T\bigl(U_u(t)-U_v(t)\bigr).
		\]
		Since $T$ is invertible, $U_u(t)=U_v(t)$. Comparing the $(2,1)$ entries
		gives $iu(t,x)=iv(t,x)$, hence $u=v$. \bigskip
	\end{proof}
	
	\begin{proof}[Proof of Theorem \ref{t1}]
		For $G=\varnothing$ the reflectionless Dirac operator with spectrum $\mathbb R$ has zero potential by \cite[Theorem~6.7]{CG}, so we take $F=0$. For $G\ne\varnothing$, the frequencies in \eqref{2.13} satisfy \eqref{1.2} by \eqref{2.18}. Take $F_{q_0}=F$ as in Lemma~\ref{l10} and use the $C^m$ representative $x\mapsto F(x\alpha)$ of $q_0$. Lemma~\ref{l11} gives \eqref{1.3}.
		
		The map $(x,\mathbf t)\mapsto x\alpha+\sum_{r=2}^m t_r\beta^{(r)}$ is a continuous homomorphism from $\mathbb R^m$ to the compact abelian group $\Theta$. Finite linear combinations of the characters
		\[
		\theta\longmapsto\exp\!\left(i\sum_{j\in J_0}n_j\theta_j\right),
		\qquad J_0\subset\mathcal J\ \text{finite},\quad n_j\in\mathbb Z ,
		\]
		are uniformly dense in $C(\Theta)$ by the Stone--Weierstrass theorem. Their pullbacks are trigonometric polynomials in $(x,\mathbf t)$. Applying this to each $F_{k,\boldsymbol\ell}$ from Lemma~\ref{l11} gives joint Bohr almost periodicity of $\partial_x^k\partial_{\mathbf t}^{\boldsymbol\ell}q$ for $k+\sum_{r=2}^m r\ell_r\le m$ \cite[Chapter~VIII]{Loomis}. In particular, $\partial_x^k q$ is bounded for $0\le k\le m$. Lemma~\ref{l12} gives uniqueness for each $r$th flow and, applied successively in $t_2,\ldots,t_m$, for the common solution of \eqref{1.3}. For the cubic flow the bound on $v,v_x$ suffices by \cite[Lemma~38]{KWXZ}. \bigskip
	\end{proof}
	
	\subsection{Proof of Theorem \ref{t2}}
	
	\begin{proof}[Proof of Theorem \ref{t2}]
		Homogeneity gives positive lower linear density at every finite point of $E$. At infinity,
		\[
		|E\cap([-2R,-R]\cup[R,2R])|\ge 2R-2S_0\ge R
		\]
		for all sufficiently large $R$. Wiener's criterion therefore gives Dirichlet regularity at every boundary point, including $\infty$. Homogeneity of $E$ implies that $\Omega$ is a Widom domain
		\cite[Section~7]{EY12} and satisfies DCT \cite{SY97}. To obtain the normalization in Definition~\ref{def:dct}, apply the one-point formula \cite[(5.1)--(5.2)]{Y14} at $\xi_*$ to $(F(w)-F(z))/(w-z)$. The difference quotient belongs to $\mathcal N_+(\Omega)$ by
		\eqref{4.1}, since it is bounded near $z$ and bounded by
		$C_z(1+|F|)$ elsewhere. Moreover, the weighted integrability in the definition gives integrability of the quotient with weight $|s-\xi_*|^{-1}$.
		
		Choose $x_j\in I_j^\circ$ for each $j$. For $t\in I_k^\circ$,
		apply \cite[Lemma~3.2 and (2.6)]{JM1985} with $t$ as the point
		chosen in the $k$th gap. This gives
		\[
		\sum_{j\ne k}g_\Omega(x_j,t)\leq A_\eta,
		\]
		where $A_\eta$ depends only on the homogeneity constant $\eta$.
		
		Let $\Omega_k=\mathbb C\setminus(\mathbb R\setminus I_k^\circ)$ and $g_k=g_{\Omega_k}$. For fixed $t\in I_k^\circ$, the function
		\[
		h(z)=g_\Omega(z,t)-g_k(z,t)
		\]
		is nonnegative in $\Omega_k$ by domain monotonicity, and its singularity at $t$ is removable. On $\partial\Omega_k$ its boundary values are zero on $E$ and at most $A_\eta$ on every other gap, by the preceding estimate with an arbitrary point in that gap. At infinity $0\le h\le g_\Omega(\cdot,t)\to0$. The maximum principle gives $h\le A_\eta$. Hence for $t\in I_k^\circ\setminus\{x_k\}$,
		\[
		\sum_j g_\Omega(x_j,t)\le g_k(x_k,t)+2A_\eta .
		\]
		The single-gap identity in the proof of Lemma~\ref{l1} gives
		\[
		\frac1\pi\int_{I_k} g_k(x_k,t)\,dt
		=\sqrt{a_k^2-(x_k-c_k)^2}\le a_k .
		\]
		Tonelli's theorem now yields
		\[
		\begin{aligned}
			\sum_j(\mathcal K_\Omega1)(x_j)
			&=\frac1\pi\sum_k\int_{I_k^\circ}\sum_jg_\Omega(x_j,t)\,dt\\
			&\le\sum_k a_k+\frac{2A_\eta}{\pi}|G|
			=\left(1+\frac{4A_\eta}{\pi}\right)S_0 .
		\end{aligned}
		\]
		Since $0<\omega\le1$, the same upper bound holds for $\sum_j v_\omega(x_j)$. Taking the supremum independently over $x_j$ in each finite set of gaps and then increasing the finite set gives \eqref{1.4} with $\mathcal P_1(E)\le(1+4A_\eta/\pi)S_0$. \bigskip
	\end{proof}
	

\end{document}